\documentclass[12pt]{amsart}[12]

\usepackage{verbatim}
\usepackage{amssymb}  
\usepackage{amscd}
\usepackage{amsfonts}
\usepackage{amsmath}
\usepackage{amsthm}
\usepackage{times}
\usepackage{graphicx,subfigure}
\usepackage{xypic}

\newenvironment{pf*}[1]{\proof[#1]}{\endproof}
\newtheorem{Theorem}[equation]{Theorem}
\newtheorem{Corollary}[equation]{Corollary}
\newtheorem{Lemma}[equation]{Lemma}
\newtheorem{Proposition}[equation]{Proposition}

\theoremstyle{definition}
\newtheorem{Definition}[equation]{Definition}
\newtheorem{Example}[equation]{Example}

\newtheorem{Conjecture}[equation]{Conjecture}

\theoremstyle{remark}

\newtheorem{Remark}[equation]{Remark}

\numberwithin{equation}{section}
\numberwithin{figure}{section}

\newcommand{\C}{{\mathbb C}}

\newcommand{\G}{{\mathcal G}}

\newcommand{\Q}{{\mathbb Q}}

\newcommand{\N}{{\mathbb N}}

\newcommand{\LL}{\Lambda}
\newcommand{\mo}{\mathcal{M}}

\newcommand{\chains}{\mathcal{C}}
\newcommand{\LA}{\Lambda^*}

\newcommand{\vGam}{\varGamma}
\newcommand{\nc}{\newcommand}

\nc\repeps{\tau}

\begin{document}

\title{Supersolvable lattices of $J$-classes}

\author{Mah\.{i}r B\.{i}len Can}

\date{October 21, 2009}

	\begin{abstract}
	The purpose of this article is to investigate the combinatorial 
	properties of the cross section lattice of a $J$-irreducible monoid
	associated with a semisimple algebraic group of one of the types $A_n$, $B_n$,
	or $C_n$. 
	Our main tool is a theorem of Putcha and Renner which identifies the 
	cross section lattice in the Boolean lattice of subsets of the nodes of a 
	Dynkin diagram. 
	We determine the join irreducibles of the cross section lattice. 
	Exploiting this we find characterizations of the relatively complemented intervals. 
	By a result of Putcha, this determines the M\"{o}bius function for 
	$\Lambda$. We show that an interval of the cross section lattice is Boolean
	if and only if it is relatively complemented if and only if it is atomic. 
	We characterize distributive cross section lattices, showing that 
	they are products of chains. We determine which cross section lattices are supersolvable, 
	and furthermore, we compute the characteristic polynomials of these 
	supersolvable cross section lattices. 
	\end{abstract}

\maketitle

\section{\textbf{Introduction}}

	Let $K$ be an algebraically closed field, and let $G_0$ be a semisimple 
	linear algebraic group over $K$. 
	Let $\rho : G_0 \rightarrow GL(V)$ be an irreducible representation of $G_0$, 
	and let $G=K^*\cdot \rho(G_0)$ be the image of $G_0$ under $\rho$, 
	adjoined by the scalar matrices. The $J-irreducible\ monoid$ $M$ associated 
	with the pair $(G_0,\rho)$ is the Zariski closure of the group $G$ in $End(V)$ (see \cite{Renner85}).

	The group $G\times G$  acts on the monoid $M$ by 
	$(g,h)\cdot x = g x h^{-1}$, $x\in M,\ (g,h)\in G\times G$.  
	Let $J_x= GxG$ denote an orbit for some $x\in M$. 
	There is a natural partial ordering on the set $\{J_x:\ x\in M\}$ of orbits: 
	\begin{equation}
	J_e \leq J_f \iff J_e \subseteq \overline{J_f},\ e,f\in M.
	\end{equation}
	Here $\overline{J_f}$ denotes the Zariski closure of $J_f$ in $M$.

	In this manuscript we are concerned with the combinatorics of the set 
	$\{J_x:\ x\in M\}$ of orbits. 
	The above construction hints just a special case of a much more general 
	theory of algebraic monoids which has been developed during the past 30 
	years largely by M.S. Putcha and L.E. Renner. The reader can find more 
	about algebraic monoids in the books \cite{Putcha88} and \cite{Renner04}.

	In \cite{Putcha83}, Putcha shows that there exists a finite 
	{\em cross section lattice of idempotents,} denoted by $\Lambda$, for the set 
	of $G\times G$-orbits.  
	More precisely, there exists a finite set $\Lambda=\{e\in M:\ e^2=e\}$ of idempotents such 
	that a) $\bigcup \{J_e:\ e\in \Lambda\} = \bigcup \{J_x:\ x\in M\}$, b) if $x\in M$, $|J_x \cap \Lambda|=1$, c) 
	$\Lambda$ is a graded lattice with respect to the partial ordering 
 	\begin{equation}\label{Jclasses}
 	J_e \leq J_f \iff J_e \subseteq \overline{J_f} \iff e=ef=fe,
 	\end{equation}
	 with the rank function given by $rk(e) = \dim GeG$.
	It turns out (see Chapter 6, \cite{Putcha88}) 
	that one can recover basic subgroups of the group $G$ from 
	$\Lambda$: the centralizer of $\Lambda$ in $G$ is a maximal torus 
	$$T=C_G(\Lambda)= \{ g\in G:\ g e=e g,\ e\in \Lambda \},$$ 
	and the {\em right centralizer} 
	$$B=C_G^r(\Lambda)=\{g\in G: \ ge=ege\ \text{for\ all}\ e\in \Lambda\}$$ 
	of $\Lambda$ in $G$ is a Borel subgroup containing the maximal torus $T$.

	 A remarkable theorem of Putcha and Renner identifies $\Lambda$ by a 
	 subposet of the Boolean lattice of all subsets of a set $\Delta$ of simple 
	 roots for $G$. Before we state the Theorem, we explain briefly the terminology 
	 about simple roots and Coxeter diagrams (disjoint unions of Dynkin diagrams). 
	 The details can be found in \cite{Humphreys}.

\begin{figure}[h]
\begin{center}
\includegraphics[width=2.4in]{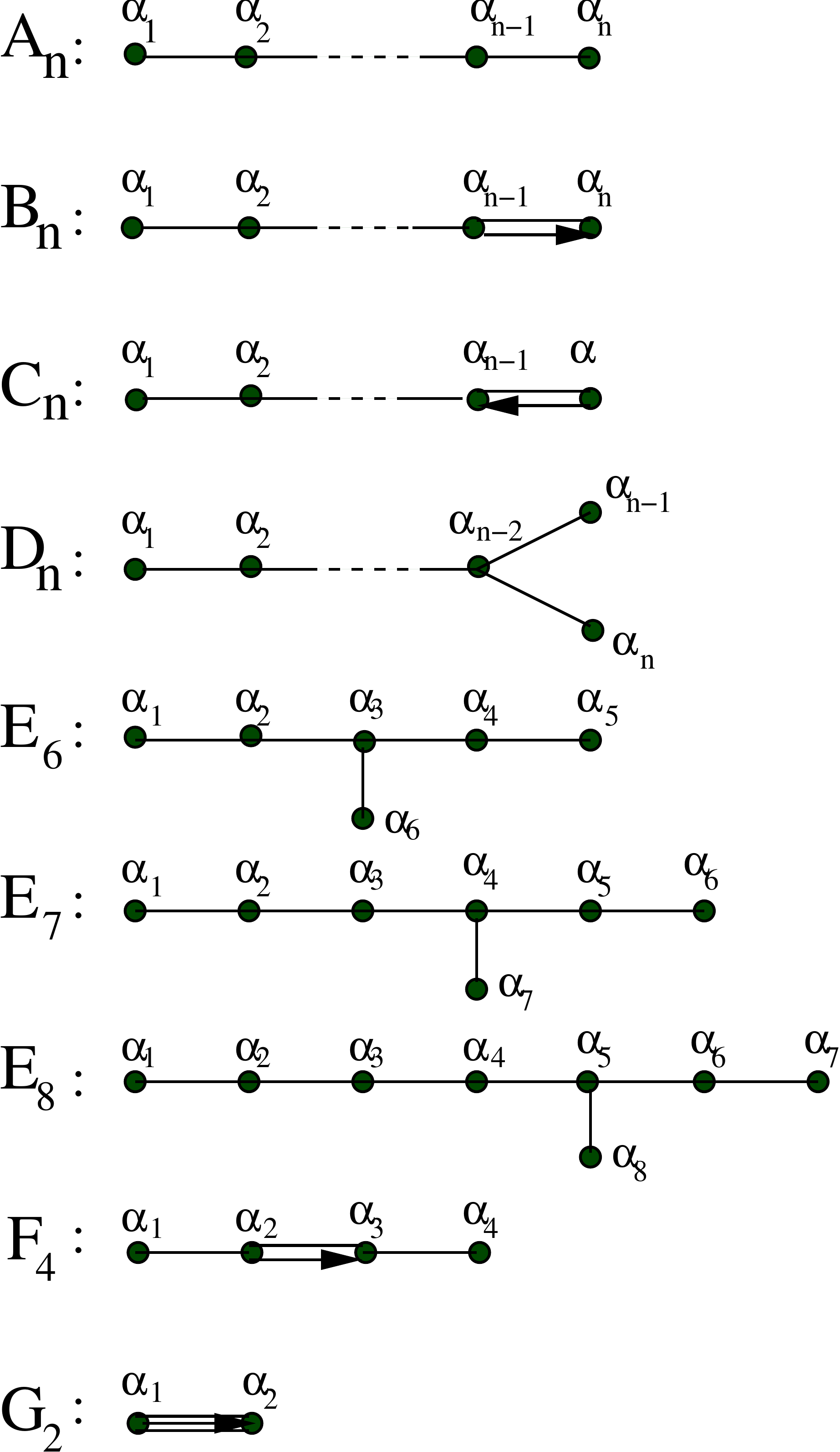}
\caption{Dynkin diagrams of the simple groups.}
\label{Dynkin}
\end{center}
\end{figure}

	Recall that simple algebraic groups are classified according to the discrete data 
	of the ``Dynkin diagrams.'' There are four infinite families of simple groups, denoted 
	by $A_n,B_n,C_n,D_n$, and five exceptional groups, denoted by $E_6,E_7,E_8,F_4,G_2$. 
	For each group in this list the {\em Dynkin diagram} is a connected ``graph'' as 
	shown in the Figure \ref{Dynkin}. Let $T$ be a maximal torus in the simple group 
	$G_0$, and let $X(T)=\text{Hom}(T,K^*)$ be the group of group homomorphisms from 
	$T$ into $K^*$. The nodes of a Dynkin diagram correspond to some special vectors, called 
	{\em simple roots}, in the vector space $V=X(T)\otimes \Q$. The set of simple roots is denoted 
	by $\Delta$. Each simple root $\alpha \in \Delta$ gives a reflection (with respect to the hyperplane 
	perpendicular to $\alpha$) on $V$, which is denoted by $\sigma_{\alpha}$. 
	The group generated by these involutions is called the 
	{\em standard reflection representation} of the Weyl group of $G$ relative to $T$. 
	Here, a  {\em Weyl group}\footnote{More generally, for an algebraic group $G$ a 
	Weyl group is defined to be the quotient group $W=N_G(T)/Z_G(T)$, where 
	$Z_G(T)=\{g\in G:\ gt=tg,\ \text{for all}\ t\in T\}$ is the centralizer of $T$ in $G$. } 
	is defined to be the quotient of  the normalizer $N_{G_0}(T)=\{g\in G:\ gT=Tg\}$by $T$.

	When the group is not simple, one still has a diagram (as well as $\Delta$), 
	called the {\em Coxeter diagram}, which is a disjoint union of the Dynkin 
	diagrams of the simple components. A subset $Y \subseteq \Delta$ is called 
	{\em connected} if it is a connected subset of the underlying graph of the Coxeter diagram.

	We turn back to the monoids. If $M$ is a $J-$irreducible monoid, 
	then $\Lambda$ has a unique minimal nonzero element, denoted 
	by $e_0$ (Lemma 7.8, \cite{Renner04}).
	\begin{Theorem} \label{T:PR88} (Putcha-Renner,  \cite{PR88})
	Let $M$ be a $J$-irreducible monoid with a cross section lattice 
	$\Lambda$. Let $T=C_G(\Lambda)$, $B=C_G^r(\Lambda)$. 
	Let $\Delta$ be the set of simple roots of $T$ relative to $B$. 
	Let $2^{\Delta}$ be the Boolean lattice of all subsets of $\Delta$. 
	Define
	$$
	\phi :  \Lambda \longrightarrow 2^{\Delta} 
	$$
	$$
	\phi(e) =  \{\alpha \in \Delta|\ \sigma_{\alpha} e = e \sigma_{\alpha} \neq e \}.
	$$
	Here, $\sigma_{\alpha}$ is the simple reflection associated with the root $\alpha$. 
	Then
	\begin{enumerate}
	\item[(i)] $\phi$ is injective and order preserving.
	\item[(ii)] $I\in 2^{\Delta}$ is in the image of $\phi$ 
	iff no connected component of $I$ (in the Coxeter graph of $\Delta$) 
	lies entirely in 	$J_0= \{ \alpha \in \Delta| \ \sigma_{\alpha} e_0 = e_0 \sigma_{\alpha}\}$ 
	where $e_0\in \Lambda - \{0\}$ is the minimal element. 
	\end{enumerate}
	\end{Theorem}

	The above theorem, which we exploit is the bridge between a 
	cross section lattice and the combinatorics of finite sets. 
	Now we can state our main results and give a brief summary of the article. 
	In the next two sections we setup the notation and give the necessary background.

	The M\"obius function of a poset $P$ is the unique function 
	$\mu: P\times P \rightarrow \N$ satisfying 
	\begin{enumerate}
	\item $\mu(x,x) =1$ for every $x\in P$,
	\item $\mu(x,y)=0$ whenever $x \nleq y$,
	\item $\mu(x,y) = - \sum_{x\leq z < y} \mu(x,z)$ for all $x<y$ in $P$.
	\end{enumerate}
	It is well known that the set of all chains of a poset is a simplicial complex 
	and the M\"obius function is the reduced Euler characteristic of the associated 
	topological space (see \cite{Stan97}). The M\"obius function
	of a cross section lattice is computed by Putcha:
	\begin{Theorem}(Putcha, \cite{Putcha04})\label{T:PM}
	Let $\Lambda$ be a cross section lattice of an algebraic monoid. 
	Let $e\leq f$ be in $\Lambda$. Then,
	\begin{equation}
	\mu(e,f)=
	\begin{cases}
	(-1)^{rk(f)-rk(e)} & \text{if}\ [e,f]\ \text{is relatively complemented,}\\
	0 & otherwise.
	\end{cases}
	\end{equation}
	\end{Theorem}
	Here relatively complemented means that for every interval 
	$[x,y]\subseteq L$, and for every $z\in [x,y]$, there exists $z'\in [x,y]$ such that  
	$z \vee z' = y$ and $z\wedge z'= x$.

	In Section \ref{S:Jirreducible} we describe a method for deciding when an 
	interval $[e,f]$ in $\Lambda$ is relatively complemented.  
	We use the notation $\LA$ to denote $\Lambda - \{0\}$.

	\begin{Theorem}\label{T:relcomp=atomic}
	We use the notation of the Theorem \ref{T:PR88}. 
	In particular, let $J_0= \{ \alpha \in \Delta| \ \sigma_{\alpha} e_0 = e_0 \sigma_{\alpha}\}$ 
	where $e_0\in \LA$ is the minimal element. 
	Then, for $e\leq f$ in $\Lambda$, the interval $[e,f]$ is relatively complemented 
	if and only if there does not exist $\alpha \in J_0$ such that 
	$\alpha \in \phi(f) - \phi(e)$ and $\sigma_{\alpha} \sigma_{\beta} =  \sigma_{\beta}\sigma_{\alpha}$
	for every $\beta \in \phi(e)$. 
	\end{Theorem}

	A graded lattice $L$ is called {\em (upper) semimodular} if the following inequality 
	\begin{equation}\label{E:semi}
	rk(x)+ rk(y) \geq rk( x \wedge y) + rk(x \vee y)
	\end{equation}
	holds for all $x,y \in L$. 
	It is known that $\Lambda$ of a $J-$irreducible monoid is upper semimodular. 
	The lattice $L$ is called atomic if every element is a join of atoms (elements covering $\hat{0}$).

	Recall that (see \cite{Stan97}) a semimodular lattice is relatively complemented 
	if and only if it is atomic. As a corollary, we also prove that 
	\begin{Corollary}
	An interval of $\LA$ is relatively complemented if and only if it is atomic if 
	and only if it is isomorphic to a Boolean lattice. 
	\end{Corollary}

	Recall that a lattice $Z$ is called {\em distributive} if for every $a,b$ and $c$ 
	from $Z$, $a\wedge (b \vee c) = (a \wedge b) \vee (a \wedge c)$ hold. 
	It is observed in \cite{PR88} that $\Lambda$ is a distributive lattice whenever 
	$\Delta - J_0$ is a connected subset of the Coxeter graph of $\Delta$, 
	where $\Delta$ and $J_0$ are as before.  In Section \ref{S:flagsymm}, we prove:

	\begin{Theorem}
	Let $\Lambda$ be the cross section lattice of a $J-$irreducible monoid 
	of the pair $(G_0,\rho)$, where $G_0$ is a simple algebraic group of one 
	of the types $A_n, B_n$ or $C_n$. Let $J_0$ be as in the 
	Theorem \ref{T:PR88}, and let $\Delta$ be the set of simple roots. 
	Suppose $|J_0|>1$.  
	Then, the followings are equivalent.
	\begin{enumerate}
	\item $\LA$ is a distributive lattice,
	\item $\LA$ is isomorphic to a product of chains,
	\item $\LA$ is isomorphic to a sublattice of a Boolean lattice.
	\end{enumerate}
	\end{Theorem}

	A particularly interesting class of distributive cross section lattices are 
	studied by Renner in \cite{Renner09}. 

	Let $M$ be a $J-$irreducible monoid, and let $\Lambda$ be the cross 
	section lattice. Let $T=C_G(\Lambda)$ be the maximal torus associated with 
	$\Lambda$, and let $W = N_G(T)/T$ be the Weyl group as before. 
	Then, $W$ acts on $V=X(T)\otimes \Q$ via its reflection representation, 
	where $X(T)=\text{Hom}(T,K^*)$. Let $\mu \in V$ be a vector which is in 
	general position with respect to lines determined by the simple roots 
	$\Delta \subseteq V$. Let $\mathcal{P}_{\mu} = Conv(W\cdot \mu) \subseteq V$ 
	be the polytope obtained by taking the convex hull of the set of points 
	$\{w\cdot \mu:\ w\in W$\} in $V$. Let $P=P_{\mu}$ be the (projective) toric 
	variety associated with this polytope.

	Assume, temporarily, that $K=\C$. Let $O$ be a variety over $K$ of dimension $n$. 
	The variety $O$ is called {\em rationally smooth} at $x\in O$, if for every point 
	$y\in U\subseteq O$ in a neighborhood of $x \neq y$ the cohomology groups 
	$H_y^i(O)=H^i(O; O- \{y\}, \Q)$ vanish except at the top dimension, 
	which is 1 dimensional. In other words,  
	$$
	H^i_y(O) = 0\ \text{for}\ i\neq 2n,\ \text{and}
	$$
	$$
	H^{2n}_y(O) = \Q.
	$$
	In \cite{Renner09}, Renner calls a subset 
	$J\subset S=S(\Delta)=\{\sigma_{\alpha}:\ \alpha \in \Delta \}$ 
	of the simple generators of $W$ by {\em combinatorially smooth}, if the toric variety 
	$O=P $ associated with the polytope $\mathcal{P}_{\mu}$ is rationally smooth 
	at every point. In his paper, Renner lists all combinatorially smooth subsets for a Weyl group.  
	Looking at his table one sees that	\footnote{for the sake of space we mention 
	type $A_n$, $n\geq 2$ only.} the subset $J$ is combinatorially smooth if 
	$J\subseteq \{\alpha_1,...,\alpha_n\}$  is equal to one of the followings: 

	\begin{enumerate}
	\item $J_0 = \emptyset$,
	\item $J_0= \{\alpha_1,...,\alpha_i\}$ where $1\leq i<n$,
	\item $J_0= \{\alpha_j,...,\alpha_n\}$ where $1<j \leq n$,
	\item $J_0= \{\alpha_1,...,\alpha_i,\alpha_j,...,\alpha_n\}$ where $ 1\leq i,j \leq n$ and $j-i \geq 3$.
	\end{enumerate}

\begin{figure}[ht]
	\begin{center}
	\includegraphics[ scale=0.3]{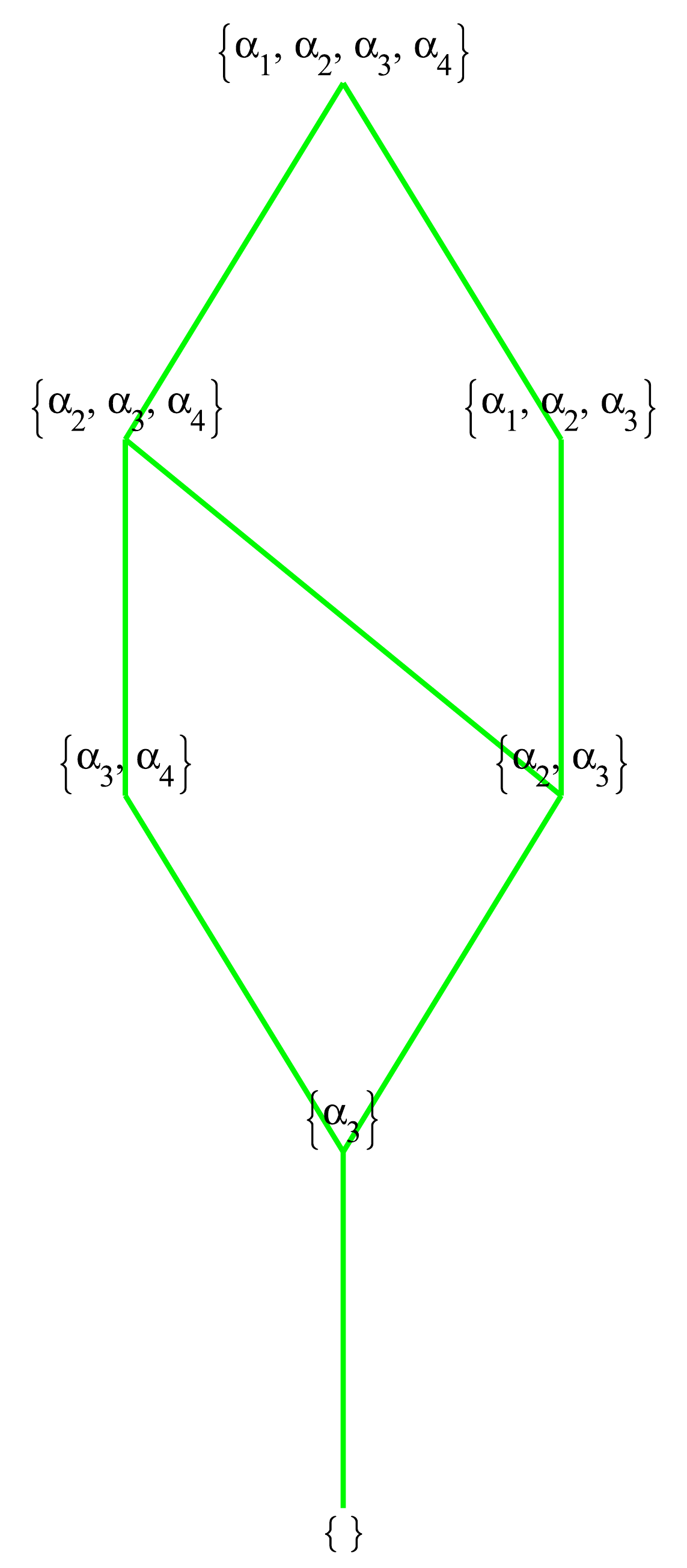}
	\caption{Combinatorially {\em non-smooth} but distributive $\LA$.}
	\label{F:nonsmooth}
	\end{center}
	\end{figure}

	In other words in type $A_n$ $P$ is a rationally smooth toric variety, then 
	$\LA$ is distributive.  However, the converse need not to be true:
	In Figure \ref{F:nonsmooth} we have an example of a cross section 
	lattice such that the associated toric variety is not rationally smooth.

	This motivates us to understand the structure of a distributive $\LA$ more precisely. 
	Thus, we make the following conjecture:

	\begin{Conjecture}\label{C:chains}
	Let $\Lambda$ be a distributive cross section lattice of the $J$-irreducible monoid of 
	the pair $(G_0,\rho)$, where $G_0$ is a simple algebraic group of one of the types 
	$A_{n-1}, B_{n-1}$ or $C_{n-1}$. Let $J_0$ and $\Delta = \{\alpha_1,...,\alpha_{n-1} \}$ 
	be the set 	of simple roots as in Theorem above. Suppose that $J_0$ has connected 
	components $J_0^{(1)}=\{\alpha_1,\alpha_2,...,\alpha_k\}$ 
	and  $J_0^{(2)}=\{\alpha_l,\alpha_{l+1},...,\alpha_{n-1}\}$ (with possibility that 
	$J_0^{(1)}=\emptyset\ \text{or that}\ J_0^{(2)}=\emptyset	$). 
	Then $\LA$ is isomorphic to the product of the chains:

	$$C_{k+2} \times C_{n-l+2} \times C_2 ^{l-k-3}  = 
	C_{k+2} \times C_{n-l+2} \times  C_2 \times \cdots \times C_2.$$
	\end{Conjecture}

	Suppose that $L$ is a finite lattice, and let $\vGam$ be a maximal 
	chain such that for every chain $\vGam'$ of $L$, the 
	sublattice generated by $\vGam$ and $\vGam'$ is distributive. 
	Such a lattice $L$ is called {\em supersolvable},
	and the maximal chain $\vGam$ is called an {\em $\mathcal{M}$-chain}.   
	This notion about lattices, introduced and used by Stanley \cite{Stan72} is a 
	direct generalization of the supersolvability
	of a group (where there is  a series of normal subgroups with all the factors are cyclic groups).
	In Section \ref{S:SS}, we determine the supersolvable cross section lattices with respect to $J_0$. 
	More precisely, we prove that

	\begin{Theorem}\label{T:ss}
	
	Let $\Delta = \{\alpha_1,...,\alpha_n\}$ be the set of simple roots and 
	$J_0$ be as before. Then, the cross section lattice $\Lambda$ is 
	supersolvable if and only if a connected component of $J_0$ is either 
	a singleton $\{\alpha_i\}$ for some $\alpha_i \in \Delta$, 
	or is a subset $\widetilde{J_0} \subseteq \Delta$ whose 
	induced subgraph in the Dynkin diagram 
	is connected and contains an end node of $\Delta$.

	Here, a simple root $\alpha \in \Delta$ is called an end-node if there exists a unique 
	$\alpha'$ in $\Delta$ such that 
	$\sigma_{\alpha} \sigma_{\alpha'} \neq \sigma_{\alpha'} \sigma_{\alpha}$ 
	where $\sigma_{\alpha'}$ is the simple reflection associated with the root $\alpha'$.

	\end{Theorem}

	We prove this theorem by exhibiting an $\mathcal{M}$-chain $\vGam$ 
	as in the definition of the supersolvability. 

	The {\em characteristic polynomial} of  a graded poset $P$ with $\hat{0}$ and $\hat{1}$ is 
	\begin{equation}
	p(x,P) = \sum_{x\in P} \mu(\hat{0},x) x^{rk(\hat{1})-rk(x)}.
	\end{equation}
	Here $rk: P\rightarrow \N$ is the rank function, and $\mu$ is the M\"obius function on $P$. 
	
	\begin{Theorem}(Stanley, \cite{Stan72}) \label{T:Stanley72}
	Let $L$ be a semimodular supersolvable lattice. 
	Let $\hat{0} = x_0 < x_1 <\cdots < x_n = \hat{1}$ be an $\mathcal{M}$-chain. Then,
	\begin{equation}
	p(x,L) = (x - a_1)(x -a_2) \cdots (x - a_n), 
	\end{equation}
	where $a_i$ is the number of atoms $u\in L$ such that $u \leq x_i$ and $u \nleq  x_{i-1}$.
	\end{Theorem}
	In Section \ref{S:CP}, using the above theorem of Stanley we prove  
	\begin{Theorem}

	Let $\LA$ be a supersolvable cross section lattice of a $J-$irreducible monoid 
	of the pair $(G_0,\rho)$, where $G_0$ is a simple algebraic group of one of the types 
	$A_n, B_n$ or $C_n$. Let $J_0\subseteq \Delta$ be as before, and let $n=|\Delta|$.
	Then the characteristic polynomial of $\LA$ is 
	\begin{equation}
	p(x,\LA)=  x^{|J_0|} (x-1)^{n-|J_0|}.
	\end{equation}
	\end{Theorem}

	Based on computer experiments, we conjecture that 
	\begin{Conjecture}
	Let $\LA$ be the cross section lattice of a ($J-$irreducible) monoid $M$. Then
	\begin{equation}
	p(x,\LA)=  x^{|J_0|} (x-1)^{n-|J_0|}.
	\end{equation}
	\end{Conjecture}

	Even though we state our theorems only for the 
	$J$-irreducible monoids of the types $A_n$, $B_n$ and $C_n$ the reader 
	will not have difficulty in proving similar results for the other types.

	\section{\textbf{Preliminaries}}

	\subsection{Poset terminology}\label{S:posetterminology}

	The unexplained terminology about posets and 
	(quasi) symmetric functions can be found in the books \cite{Stan97} 
	and \cite{Stan97II}. 
	We assume that all posets and lattices are finite, and furthermore have 
	a minimal element, denoted by $\hat{0}$, as well as a maximal element 
	which is denoted by $\hat{1}$.

	Let $\Delta$ be a finite set (which is going to stand for the set of simple roots).
	 We use the notation $2^{\Delta}$ to denote the Boolean lattice of all subsets 
	 of $\Delta$ with the partial order of inclusion.

	For any integer $m>0$ we denote $\{1,...,m\}$ by $[m]$. 
	An increasing sequence of numbers $m_1<m_2<\cdots < m_k$ is 
	denoted by $\{m_1,...,m_k\}_<$.

	A {\em composition} $\gamma$ of a positive integer $m$ is an ordered sum 
	of positive integers whose sum is $m$. 
	The set of all compositions of $m$ is denoted by $Comp(m)$. 
	A {\em partition} $\lambda$ of a positive integer $m$ is an unordered sum of 
	positive integers whose sum is $m$. 
	There is a $1-1$ correspondence between subsets of $\{1,...,m-1\}$ 
	and compositions of $m$ which is defined by sending 
	$\gamma=\gamma_1+\cdots + \gamma_k $ to 
	$S_{\gamma}=\{ \gamma_1 , \gamma_1+\gamma_2, ... , \gamma_1+\cdots \gamma_{k-1}\}$.  
	For a composition $\gamma$ of $m$, $\lambda=\lambda(\gamma)$ 
	denotes the partition of $m$ obtained by rewriting the entries of $\gamma$ 
	in a decreasing order.

	Let $P$ be a graded poset of rank $n$ with $\hat{0}$ and $\hat{1}$.  
	Equivalently, $P$ has the smallest and the largest elements and every 
	maximal chain in $P$ is of the same length $n$.  Let $rk: P \rightarrow \N$ 
	be the {\em rank function}, so that $rk(x)\in \N,\ x\in P$ is the length of the 
	maximal chain from $\hat{0}$ to $x$. If $x\leq y$ in $P$, then the length of 
	the interval $[x,y]=\{z\in P:\ x\leq z \leq y\}$ is defined to be 
	
	\begin{equation}
	rk(x,y) = rk(y)- rk(x).
	\end{equation}
	
	A graded poset $P$ of rank $n$ is called {\em locally rank symmetric} 
	if every interval of $P$ is {\em rank-symmetric}, i.e., has the same number 
	of elements of rank $i$ as of corank $i$ for all $i$.  
	Here corank of an element $x\in P$ means $n-rk(x)$.

	The set of all chains in an interval $I$ (resp. in $P=[\hat{0},\hat{1}]$) of 
	$P$ is denoted by  $\chains(I)$ (resp. by $\chains$). 
	Similarly, the set of all max chains in $I$ (resp. in $P=[\hat{0},\hat{1}]$) of 
	$P$ is denoted by  $\chains_{max}(I)$ (resp. by $\chains_{max}$).

	The {\em rank set} $rk(\tau)$ of a chain $\tau=(t_1<\ldots < t_i) \in \chains$ is 
	defined to be
	
	\begin{equation}
	rk(\tau) = \{m_i:\ m_i = rk(t_i)\}_<.
	\end{equation}
	
	Let $n$ be the rank of $P$, and let $I,J \subseteq \{1,\ldots,n-1\}$. Define
	\begin{eqnarray}\label{E:alpha}
	\alpha_P(I) &=& |\{  \tau :\  rk(\tau)= I \} |,\\
	\beta_P(J) &=& \sum_{I\subseteq J} (-1)^{|J - I|} \alpha_P(I).
	\end{eqnarray}

	The function $\alpha_P$ (from the set of all subsets of $[rk(P)]$ into $\N$) is 
	sometimes called the {\em flag $f$-vector} of $P$.

	The formal power series 
	\begin{equation}
	F_P(x) = \sum_{\hat{0}=t_0\leq t_1\leq\cdots\leq t_{k-1}<t_k
	    =\hat{1}} x_1^{rk(t_0,t_1)}x_2^{rk(t_1,t_2)}\cdots
	    x_k^{rk(t_{k-1},t_k)}, \label{1} 
	\end{equation}
	where the sum is over all multichains from $\hat{0}$ to $\hat{1}$
	such that $\hat{1}$ occurs exactly once is 
	introduced by Ehrenborg in \cite{E97}.
	In \cite{Stan96}, Stanley shows that 
	  $$ F_P= \sum_{I\subseteq [n-1]} \beta_P(I)F_{I,n},$$
	  where 
	  \begin{equation}\label{Gessel}
	   F_{I,n} = \sum_{j_1 \leq \cdots \leq j_n \atop i\in I \Rightarrow j_i < j_{i+1}} x_{j_1} \cdots x_{j_n},
	  \end{equation}
	are the {\em fundamental quasi-symmetric functions}.
	A poset $P$ is called {\em flag-symmetric} if the function $F_P$ is a symmetric function.

	Let $L$ be a lattice (hence joins ``$\vee$'' and the meets ``$\wedge$'' exist). 
	$L$ is called {\em relatively complemented}, if for every interval $[x,y]\subseteq L$, 
	and for every $z\in [x,y]$, there exists $z'\in [x,y]$ such that  $z \vee z' = y$ and $z\wedge z'= x$. 
	A lattice $L$ is called {\em atomic} if every element of $L$ is the join of atoms 
	(elements covering $\hat{0}$)  of $L$.

	A graded lattice $L$ is called {\em (upper) semimodular} if the following inequality 
	\begin{equation}\label{E:semi}
	rk(x)+ rk(y) \geq rk( x \wedge y) + rk(x \vee y)
	\end{equation}
	holds for all $x,y \in L$. Equivalently (Birkhoff's characteriztion): 
	if $x$ covers $x\wedge y$, then $x \vee y$ covers $y$. 
	In an (upper) semimodular lattice $L$ being relatively complemented is equivalent 
	to being atomic (see Proposition 3.3.3., \cite{Stan97}).

	A graded lattice $L$ is called modular if the inequality in (\ref{E:semi}) is equality for all 
	$x,y\in L$. A lattice $L$ is called distributive if for all $x,y,z\in L$ the following holds 
	\begin{equation}
	x \wedge ( y \vee z ) = (x \wedge y) \vee (x \wedge z ).
	\end{equation}
	It is easy to show that a distributive lattice is modular.

	A pair $(a,b)$ of elements of $L$ is called a {\em modular pair}, if the following holds
	\begin{equation}
	c \leq b \Rightarrow c \vee (a \wedge b ) = (c \vee a )\wedge b.
	\end{equation}

	It is convenient to denote a modular pair by $a \mo b$.  An element $b\in L$ is called 
	{\em right modular} if $x\mo b$ for every $x\in L$. Similarly, an element $a\in L$ is called 
	{\em left modular} if $a \mo x$ for every $x\in L$. 
	It is easy to see that a lattice is modular if and only if every ordered pair of elements of $L$ is modular.

	\subsection{Reductive groups and monoids}
	
	The purpose of this section is to introduce the notation of the reductive monoids. We focus mainly on 
	the orbit structure of these monoids. For more details the reader should consult \cite{Renner04} or
	\cite{Putcha88}. For an introduction, with many explicit examples,
	we suggest the survey by Solomon \cite{Sol95}. For more background on the theory of 
	algebraic groups, we suggest Humphrey's book \cite{Humphreys}.

	Our basic list of notation for groups is as follows. 
	\begin{align}\label{A:groupnotation}
	G &= \text{reductive group,}\\
	B &= \text{Borel subgroup,}\\
	T &= \text{maximal torus contained in $B$,}\\
	W &= N_G(T)/T = \text{Weyl group of}\ (G,T),
	\end{align}

	An algebraic monoid is a variety $M$ together with 
	\begin{enumerate}
	\item an associative morphism
	\begin{equation*}
	m: M \times M \rightarrow M
	\end{equation*}
	\item a unity $1\in M$ with $m(1,x)=m(x,1)=x$ for all $x\in M$. 
	\end{enumerate}

	The set $G=G(M)$ of invertible elements of $M$ is an algebraic group. If $G$ is a reductive group
	and $M$ is irreducible $M$ is called a \textit{reductive monoid}. It turns out \cite{Rittatore98}
	that the reductive monoids are exactly the affine, two-sided embeddings of connected reductive groups.

	Many interesting examples of (reductive) monoids are obtained as the 
	Zariski closures $\overline{K\cdot G}\subseteq End(V)$
	of (representations of) reductive groups in a space of endomorphisms on a linear space $V$.

	In a reductive monoid, the data of the Weyl group $W$ of the reductive group $G$ 
	and the set $E(\overline{T})$ of idempotents of the embedding $\overline{T} \hookrightarrow M$ 
	combine to become a finite inverse semigroup $R=\overline{N_G(T)}/T \cong W\cdot E(\overline{T})$ 
	with unit group $W$ and idempotent set $E(R)=E(\overline{T})$.
	The inverse semigroup $R$, controls the Bruhat decomposition
	for the algebraic monoid $M$. Recall that the Bruhat-Chevalley order on the Weyl group 
	$W$ is defined by
	\begin{equation*}
	x \leq y\ \mbox{iff}\ BxB \subseteq \overline{ByB},
	\end{equation*}
	where $B$ is a Borel subgroup of $G$. 
	Similarly, on the Renner monoid $R$ of a reductive monoid $M$, the Bruhat-Chevalley order is 
	defined by 
	\begin{equation}\label{E:BRordering}
	\sigma \leq \tau\ \mbox{iff}\ B\sigma B \subseteq \overline{B\tau B}.
	\end{equation}
	One observes that the induced poset structure on $W$, which is induced from $R$ 
	is the same as the original Bruhat poset structure on $W$.

	Let $T$ be a maximal torus and $E(\overline{T})$ be the set of idempotent elements 
	in the (Zariski) closure of $T\subseteq G$ in the monoid $M$. Similarly, let us denote by 
	$E(M)$ the set of idempotents in the monoid $M$. 
	Plainly we have $E(\overline{T}) \subseteq E(M)$. There is a canonical partial
	order $\leq$ on $E(M)$ (hence on $E(\overline{T})$) defined by  
	\begin{equation}\label{E:crossorder}
	e\leq f\ \Leftrightarrow   \ ef=e=fe.
	\end{equation}
	Notice that $E(\overline{T})$ is invariant under the conjugation action of the Weyl group $W$. 
	We call a subset $\Lambda \subseteq E(\overline{T})$ as a \textit{cross-section lattice} if $\LL$ 
	is a set of representatives for the $W$-orbits on $E(\overline{T})$ and the bijection 
	$\LL \rightarrow G \backslash M / G$ defined by 
	$e \mapsto GeG$ is order preserving. It turns out that we can write 
	$\LL=\LL(B) = \{ e\in E(\overline{T}):\ Be=eBe\}$
	for some unique Borel subgroup $B$ containing $T$.  
	The partial order given by (\ref{E:crossorder}) on $E(\overline{T})$ (hence on $\Lambda$) 
	agrees with the Bruhat-Chevalley order (\ref{E:BRordering}) on the Renner monoid.

	The decomposition $M= \bigsqcup_{e\in \LL} GeG$,
	of a reductive monoid $M$ into its $G\times G$ orbits, has a counterpart on the Renner
	monoid $R$ of $M$. Namely, the finite monoid $R$ can be written as a disjoint union 
	\begin{equation}\label{E:RennerDecomposition}
	R = \bigsqcup_{e\in \LL} WeW
	\end{equation}
	of $W\times W$ orbits, parametrized by the cross section lattice.

	It is shown by Putcha \cite{Putcha01} that each orbit $WeW$, for $e\in \LL$, 
	is a lexicographically shellable poset. Notice, as a special case, that if 
	$e\in \Lambda$ is the identity element of $G$, then the orbit $WeW$ is the Weyl
	group $W$, and the lexicographically shellability of the Coxeter
	groups is well known \cite{BjWs82}, \cite{Proc82}. 
	The author shows in \cite{Can08} that the Renner monoid (so called {\em rook monoid})
	of the monoid of $n\times n$ matrices is lexicographically shellable. 
	The question of shellability of a Renner monoid in general is still an open problem.

	It is known that $E(\overline{T})$ is a relatively complemented lattice, anti-isomorphic 
	to a face lattice of a convex polytope.  Let $\Lambda$ be a cross section lattice in $E(\overline{T})$.  
	The Weyl group of $T$ (relative to $B= C_G^r(\Lambda)$) acts on $E(\overline{T})$, 
	and furthermore
	\begin{equation}\label{E:ELAMB}
	E(\overline{T})= \bigsqcup_{w\in W} w \Lambda w^{-1}.
	\end{equation}

	 It is shown by Putcha [Corollary 8.12, \cite{Putcha88}] that if 
	 $\emptyset \neq \vGam \subseteq E(\overline{T})$ is such that $(\vGam, \leq)$ 
	 is a relatively complemented lattice with all maximal chains having length equal to 
	 $\dim T$, then $\vGam = E(\overline{T})$. Therefore, the cross section lattice 
	 $\Lambda \subseteq E(\overline{T})$ is relatively complemented if and only if 
	 $\Lambda = E(\overline{T})$, and this is possible if and only if $W$ is trivial.

	A reductive monoid with $0$ is called $J-$irreducible if there exists a unique minimal, 
	nonzero $G\times G$ orbit (= a $J$-class). This explains the terminology.  
	Notice that this is equivalent to the fact that $\Lambda$ has a unique minimal nonzero 
	element $e_0 \in \Lambda$.  One of the many reasons to be interested in $J-$irreducible 
	monoids is a theorem of Renner \cite{Renner85}, which shows that a monoid $M$ is 
	$J-$irreducible if and only if there exists a rational representation $\rho: M \rightarrow End(K^n)$ 
	such that $\rho$ is a finite morphism, and $K^n$ is an irreducible module for $M$.


	\section{\textbf{The cross section of a $J-$irreducible monoid}}\label{S:Jirreducible}

	Let $\alpha \in \Delta$ be a simple root. 
	We denote by $\sigma_{\alpha}$ the simple reflection corresponding to $\alpha$. 
	Let $\Lambda^*=\Lambda - \{0\}$.

	The following theorem of  Putcha and Renner is one of the crucial steps in determining 
	the structure of the cross section lattice of a $J-$irreducible monoid.

	\begin{Theorem} (Putcha and Renner, Theorem 4.13, \cite{PR88})\label{T:PR4.13}
	We use the notation of the Theorem \ref{T:PR88}. Let $e\in \LA$. 
	Then there is a one to one correspondence between 
	$\{\alpha \in \Delta:\ \sigma_{\alpha} e \neq e \sigma_{\alpha} \}$ in $\Delta$ 
	and $\{f\in \Lambda:\ f\ \text{covers}\ e \}$. Furthermore, $\alpha$ corresponds to the 
	unique $f=f_{\alpha}$ with $\phi(f)=\phi(e) \cup \{\alpha\}$. 
	\end{Theorem}

	By Theorem \ref{T:PR88}, we see that $\phi(e \wedge f)$ is the largest subset $U$ 
	of $\phi(e) \cap \phi(f)$ such that there exists an idempotent $h\in \Lambda$ with 
	$U=\phi(h)$. Therefore, the following corollary of 
	Theorem \ref{T:PR4.13} (and Theorem \ref{T:PR88}) is straightforward.

	\begin{Corollary}\label{C:PR88}
	
	Let $e,f \in \Lambda$, and let $U=\phi(e),\ V=\phi(f)$. Then, 
	\begin{enumerate} 
	\item $\alpha \in \phi(e\wedge f)$ if and only if $\alpha \in U\cap V$, and $J_0$ 
	does not contain the connected component in $U\cap V$ of $\alpha$. 
	\item If  $V\cap J_0 =\emptyset$, then $\phi(e\wedge f)=U\cap V$.  
	\item If  $V\cap J_0 =\emptyset$, then for any subset $Y \subseteq V$ 
	there exists $h\in \Lambda$ such that $h \leq f$ and $Y =\phi(h)$. 
	\end{enumerate}
	
	\end{Corollary}

	The following, which was already known to the experts (see Remark 3.6, \cite{Putcha01}) 
	is also a straightforward corollary of the above considerations.

	\begin{Corollary} \label{C:rank}
	Let $e_0 \in \Lambda$ be the minimal nonzero idempotent. 
	Then, $\Lambda$ is an (upper) semimodular lattice. Furthermore, 
	the value of the rank function $e \mapsto \dim (eT)$ on $\Lambda$ is equal to 
	$|\phi(e)|=| \{ \alpha \in \Delta:\ \sigma_{\alpha} e = e \sigma_{\alpha} \neq e\}|$. 
	\end{Corollary}

	\begin{proof}
	By Theorems \ref{T:PR4.13} and \ref{T:PR88}, $\phi(e \vee f) = \phi(e) \cup \phi(f)$.  
	If $e$ covers $e\vee f$, then there exists a unique $\beta \in \phi(e) - \phi(f)$ 
	such that  $\phi(e) = \phi(e\wedge f) \cup \{\beta\}$. 
	Then, $\phi(e\vee f) = \phi(e) \cup \phi(f) = \phi(f) \cup \{\beta\}$. 
	In other words, $e\vee f$ covers $f$. 
	This proves the semimodularity. 
	The second assertion is straightforward from Theorems \ref{T:PR4.13} and \ref{T:PR88}, also. 
	\end{proof}

	The set of join irreducibles of $\LA$ is found as follows.

	\begin{Definition}\label{D:connected}
	
	Let $J \subseteq \Delta$ be a subset of the nodes of the Dynkin diagram. 
	A subset $A=\{\alpha_1,...,\alpha_k\}$ of $J$ is called  
	{\em a connected subset} if the subgraph of the Dynkin 
	diagram induced by $A$ is connected.
	Therefore, if $|A| > 1$, then for each $\alpha_i \in A$ there exists 
	$\alpha_j \neq \alpha_i$ in $A$ such that 
	$\sigma_{\alpha_i} \sigma_{\alpha_{j}} \neq \sigma_{\alpha_j} 	\sigma_{\alpha_i}$.
	A subset $A \subseteq J$ is called {\em a connected component} of $J$ if its 
	a maximal connected subset of $J$. 
	Without loss of generality we assume that $\emptyset$ is a connected subset of $J$. 
	
	\end{Definition}

	\begin{Proposition}
	Let $J_0 \subseteq \Delta$. Then, $e\in \LA$ is a join irreducible if one of the 
	following statements holds:
	\begin{enumerate}
	\item[a)] $\phi(e)$ is a singleton $\{ \beta \}$ where $\beta \in \Delta - J_0$, 
	\item[b)] $\phi(e)=A \cup \{\beta\}$, where $\beta \in \Delta - J_0$ and $A$ is 
	a connected subset of $J_0$  such that  there exists an $\alpha \in A$ with 
	$\sigma_{\alpha} \sigma_{\beta} \neq \sigma_{\beta} \sigma_{\alpha}$.
	\end{enumerate}
	\end{Proposition}

	\begin{proof}

	Let $A$ be a connected subset of $J_0$ and let $\beta \in \Delta - J_0$. 
	Since $\beta$ is an atom, It is a join irreducible. Suppose that there exist $\alpha \in A$ 
	with $\sigma_{\alpha} \sigma_{\beta} \neq \sigma_{\beta} \sigma_{\alpha}$. 
	Therefore, $A \cup \{\alpha \}$ is of the form $\phi(e)$ for some $e\in \LA$. 
	Assume for a second there exist idempotents such $e_1,e_2 \in \LA$ such that 
	$e=e_1 \vee e_2$, and neither $e_1 \leq e_2$, nor $e_2 \leq e_1$ is true. 
	Since $A\subseteq J_0$, and since $A\cup \{\alpha\} = \phi(e_1) \cup \phi(e_2)$,
	the assumption forces either $\phi(e_1)$ or $\phi(e_2)$ to have a connected 
	component lying in $A\subseteq J_0$. 
	Therefore, such $e_1$ and $e_2$ cannot exist, hence, $A\cup \{\alpha\}$ is a join irreducible.

	Next we show that these are all of the join irreducibles. It is convenient to identify 
	$\LA$ by its image in $2^{\Delta}$. 
	Obviously, any element of  $\LA$ is of the form $I=A \cup B$ for some connected 
	subset $A \subseteq J_0$ (possibly empty) and $B\subseteq \Delta$. 
	Assume that $I$ is a join-irreducible. Suppose that $A$ is maximal with respect 
	to the property that $A\subseteq I$.  Then, there must exists 
	$\beta \in (\Delta - J_0) \cap B$ such that
	$\sigma_{\beta} \sigma_{\alpha} \neq \sigma_{\alpha} \sigma_{\beta}$ for some 
	$\alpha \in A$. Since $A$ is maximal, if $I \neq A \cup \{\beta\}$, then $I- A$  
	is an element of $\LA$, so is $A \cup \{\beta\}$. 
	Then, $(I- A) \vee (A\cup \{\beta\}) = (I- A) \cup (A\cup \{\beta\}) =I$. 
	But this contradicts with our assumption that $I$ is a join-irreducible. 

	\end{proof}

	Recall that a lattice is called relatively complemented, if for every interval $[x,y]$ 
	and for every $z\in [x,y]$, there exists $z'\in [x,y]$ such that $z \vee z' = y$ and 
	$z \wedge z' = x$.  Recall also that a lattice is called atomic, if every element of the 
	lattice is a join of atoms. In a semimodular lattice these two concepts 
	are the same (see \cite{Stan97}). 
 
	Therefore, an atomic interval in $\LA$ can be characterized as in the Theorem 
	\ref{T:relcomp=atomic}.

	\noindent
	{\em Proof of Theorem \ref{T:relcomp=atomic}.}
	($\Rightarrow$) Assume contrary that there exists $\alpha \in J_0$ such that 
	$\alpha \in \phi(f) - \phi(e)$ and 
	$\sigma_{\alpha} \sigma_{\beta} =  \sigma_{\beta}\sigma_{\alpha}$
	for some $\beta \in \phi(e)$.  
	Then, by Theorem \ref{T:PR4.13}  $\phi(e) \cup \{\alpha \}$ cannot be in the 
	image $\phi(\LA)$. 
	This forces $|\phi(f) - (\phi(e) \cup \{\alpha\})|\geq 1$. 
	Let $\mathcal{A}_f$ be the set of all subsets 
	$U\subseteq \phi(f) -(\phi(e) \cup \{\alpha \})$ such that 
	$U \cup \phi(e) \cup \{\alpha \}$ lies in the interval $[\phi(e), \phi(f)]$. 
	Clearly $\mathcal{A}_f$ is non empty. Let $U_0$ be a minimal element of $\mathcal{A}_f$ 
	(with respect to inclusion ordering).  
	Clearly $U_0 \cup \phi(e) \cup \{ \alpha \}$ is a join irreducible. Since 
	$U_0 \neq \emptyset$, $U_0 \cup \phi(e) \cup \{ \alpha \}$ is not an atom. 
	Since in a relatively complemented interval a join-irreducible has to be an 
	atom, we obtain a contradiction.

\begin{figure}[h]\label{F:dortresim}
	  \centering
	  \subfigure[]
	  {
	      \includegraphics[width=1.5in]{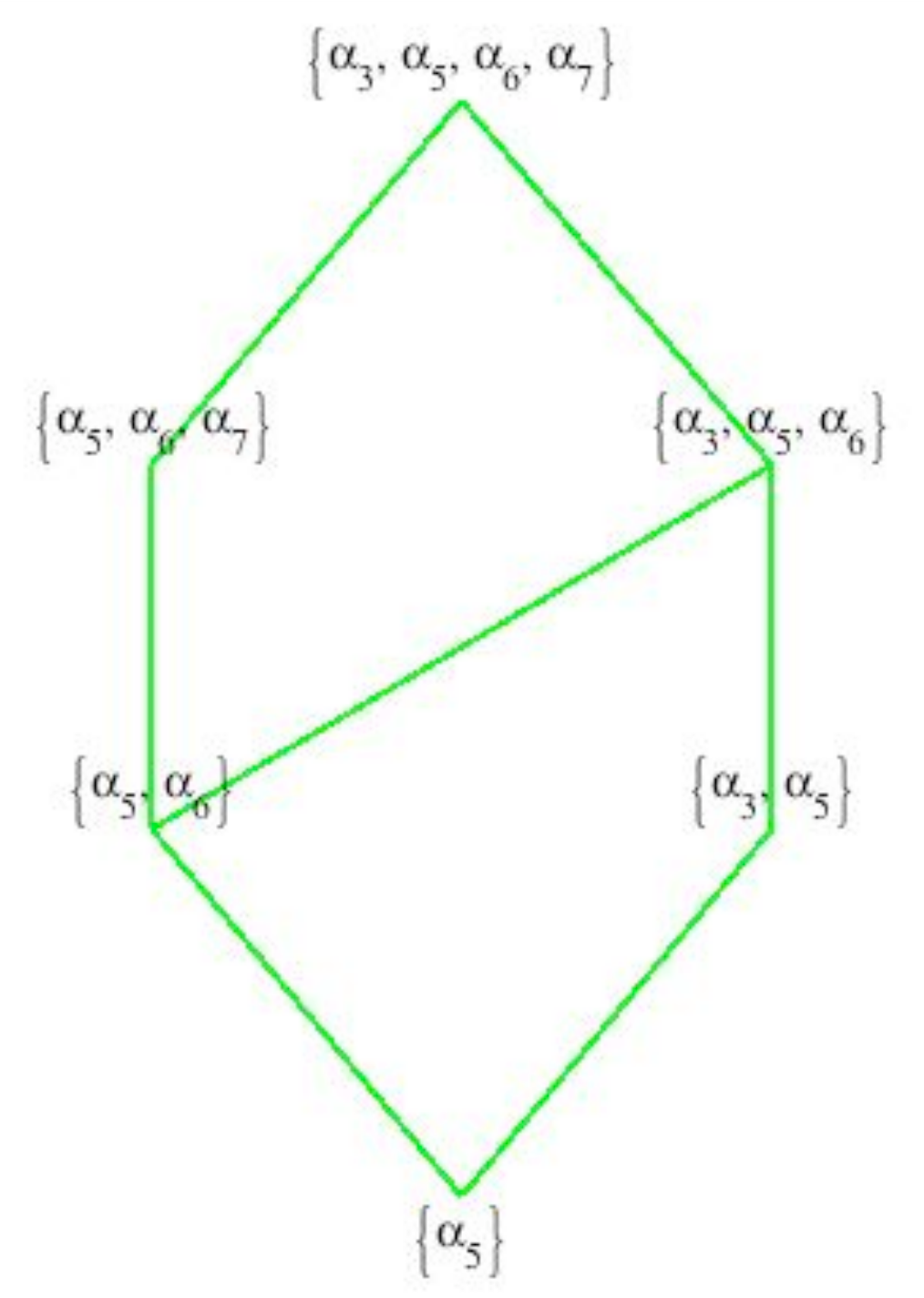}
	      \label{F125}
	  }
	   \subfigure[]
	  {
	      \includegraphics[width=2.1in]{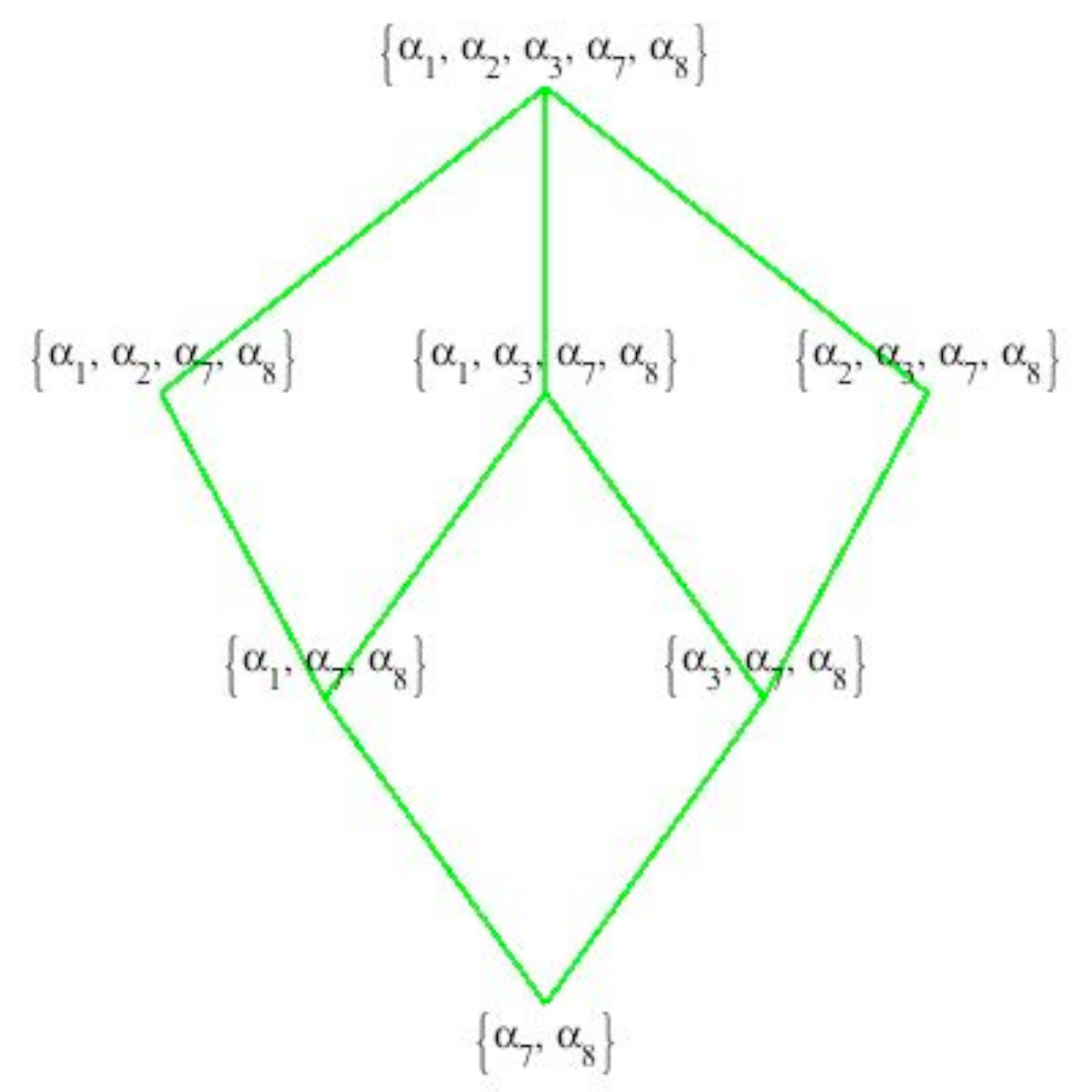}
	      \label{F123}
	  }
	     \caption{Some non-atomic subintervals.}
	  \label{F:non-atomic}
	\end{figure}

	 \noindent $Example.$ Let $\Delta = \{\alpha_1,...,\alpha_8\}$ be a set of simple roots with 
	 $\alpha_i \alpha_{i+1} \neq \alpha_{i+1} \alpha_i$, $i=1,...,7$, for $G=GL_n$. 
	 Suppose that $J_0= \{\alpha_3, \alpha_6,\alpha_7\}$. Let $\phi(e)= \{\alpha_5\}$, and let 
	 $\phi(f)=\{\alpha_3,\alpha_5,\alpha_6,\alpha_7\}$. Then, the pair $(\alpha,\alpha')$ can be 
	 chosen to be $(\alpha_7,\alpha_6)$. 
	Clearly $\phi(e) \cup \{\alpha_6,\alpha_7\} = \{\alpha_5,\alpha_6,\alpha_7\}$ is a join irreducible  
	in $[\phi(e),\phi(f)]$ but not an atom.
	Next, let $\phi(e)= \{\alpha_7,\alpha_8\}$, and let $\phi(f)=\{\alpha_1,\alpha_2,\alpha_3,\alpha_7,\alpha_8\}$. 
	Then, the pair $(\alpha,\alpha')$ can be chosen to be $(\alpha_2,\alpha_3)$. 
	Clearly $\phi(e) \cup \{\alpha_1,\alpha_2\} = \{\alpha_1,\alpha_2,\alpha_7,\alpha_8\}$ 
	is a join irreducible  in $[\phi(e),\phi(f)]$ but not an atom. 
	We depict the Hasse diagrams of these intervals in Figure \ref{F:non-atomic}.

	($\Leftarrow$) We prove the converse statement. Let $\phi(e) \subseteq 
	A \subseteq \phi(f)$ be an arbitrary subset. Let $\beta \in A - \phi(e)$ 
	be a simple root. Hence, by hypotheses, either $\beta \notin J_0$ or $\beta \in J_0$ 
	and $\sigma_{\alpha'} \sigma_{\beta} = \sigma_{\beta} \sigma_{\alpha'}$ for some 
	$\alpha' \in \phi(e)$.  In both of these cases, we see that $\phi(e) \cup \{\beta\} \in \phi(\LA)$. 
	But we can repeat this for all elements of $A- \phi(e)$. Therefore, $A \in \phi(\LA)$. 
	This shows that the interval $[e,f]$ is isomorphic to the Boolean lattice 
	$2^{\phi(f) - \phi(e)}$ which is relatively complemented, of course. 
	This finishes the proof.

	The second part of the proof the Theorem shows that

	\begin{Corollary}
	An interval of $\LA$ is relatively complemented if and only if atomic if and only if 
	isomorphic to a Boolean lattice. 
	\end{Corollary}

	As an application we rewrite the Theorem \ref{T:PM} 
	\begin{Corollary}
	Let $\Lambda$ be a cross section lattice of an $J-$irreducible monoid. Let $e\leq f$ be in $\Lambda$. Then,
	\begin{equation}
	\mu(e,f)=
	\begin{cases}
	(-1)^{rk(f)-rk(e)} & \text{if}\ [e,f]\ \text{is a Boolean lattice},\\
	0 & otherwise.
	\end{cases}
	\end{equation}
	\end{Corollary}

	\section{\textbf{Flag symmetric cross section lattices}}\label{S:flagsymm}

	Recall that a poset is called flag-symmetric if and only if the flag-quasi 
	symmetric function $F_P$ is a symmetric function. 
	In this section we answer (at least, partially) the following question:
	\vskip.5cm

	\begin{center}
	{\em Which $J$-irreducible monoids do have their cross section lattice flag-symmetric? } 
	 \end{center}
 
	 \vskip.1in

	The following observation of Stanley is useful for deciding  when a distributive lattice is flag-symmetric.

	\begin{Theorem}(Stanley, \cite{Stan96})
	Let $L$ be a finite distributive lattice. The following four conditions are equivalent.
	\begin{enumerate}
	\item $L$ is locally self dual.
	\item $L$ is locally rank-symmetric.
	\item $L$ is flag-symmetric.
	\item $L$ is a product of chains.
	\end{enumerate}
	\end{Theorem}

\begin{figure}[htp]\label{F:dortresim}
	  \centering
	  \subfigure[$J_0=\{\alpha_1,\alpha_2,\alpha_5\}$]
	  {
	      \includegraphics[width=2.3in]{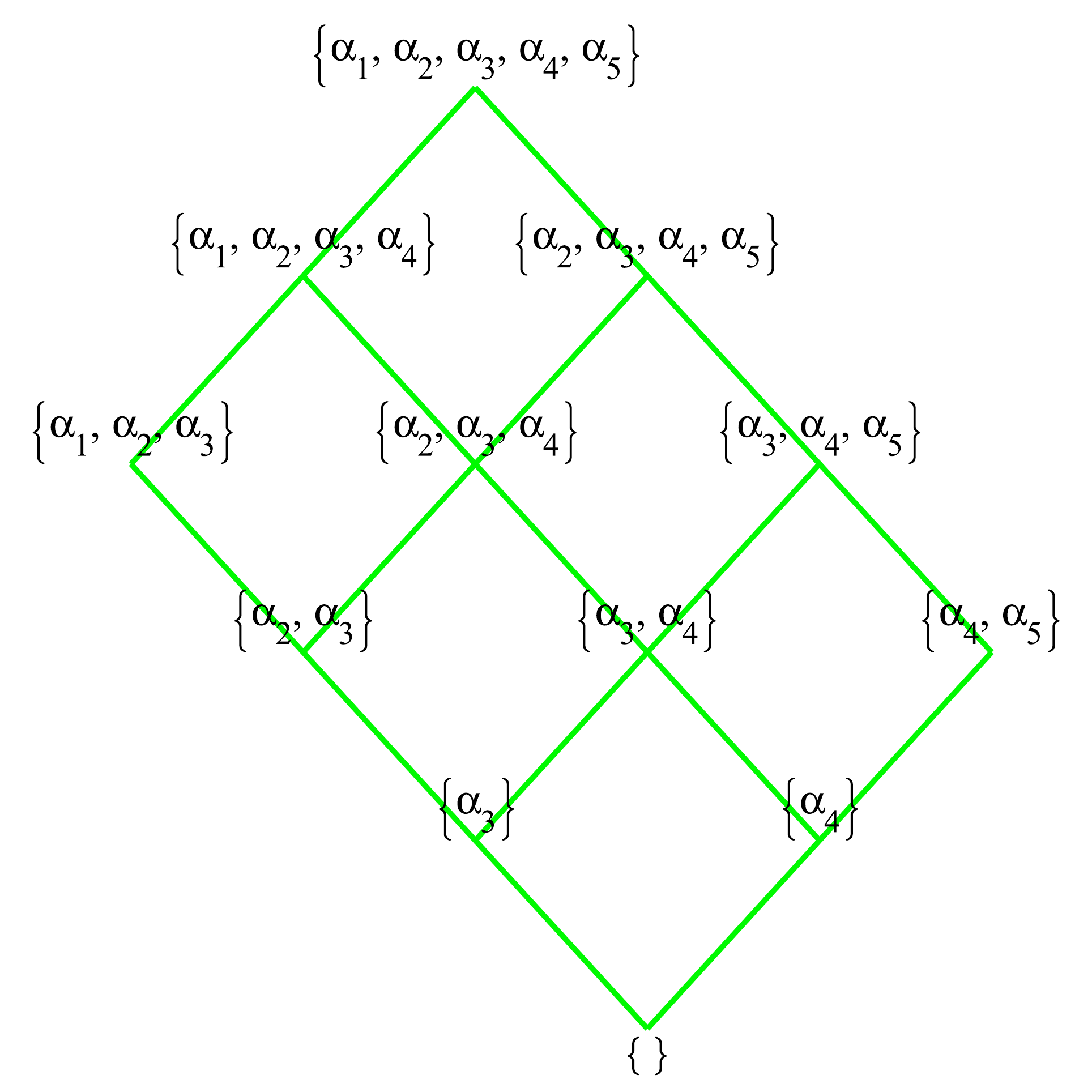}
	      \label{F125}
	  }
	   \subfigure[$J_0=\{\alpha_1,\alpha_2,\alpha_3\} $]
	  {
	      \includegraphics[width=1.3in]{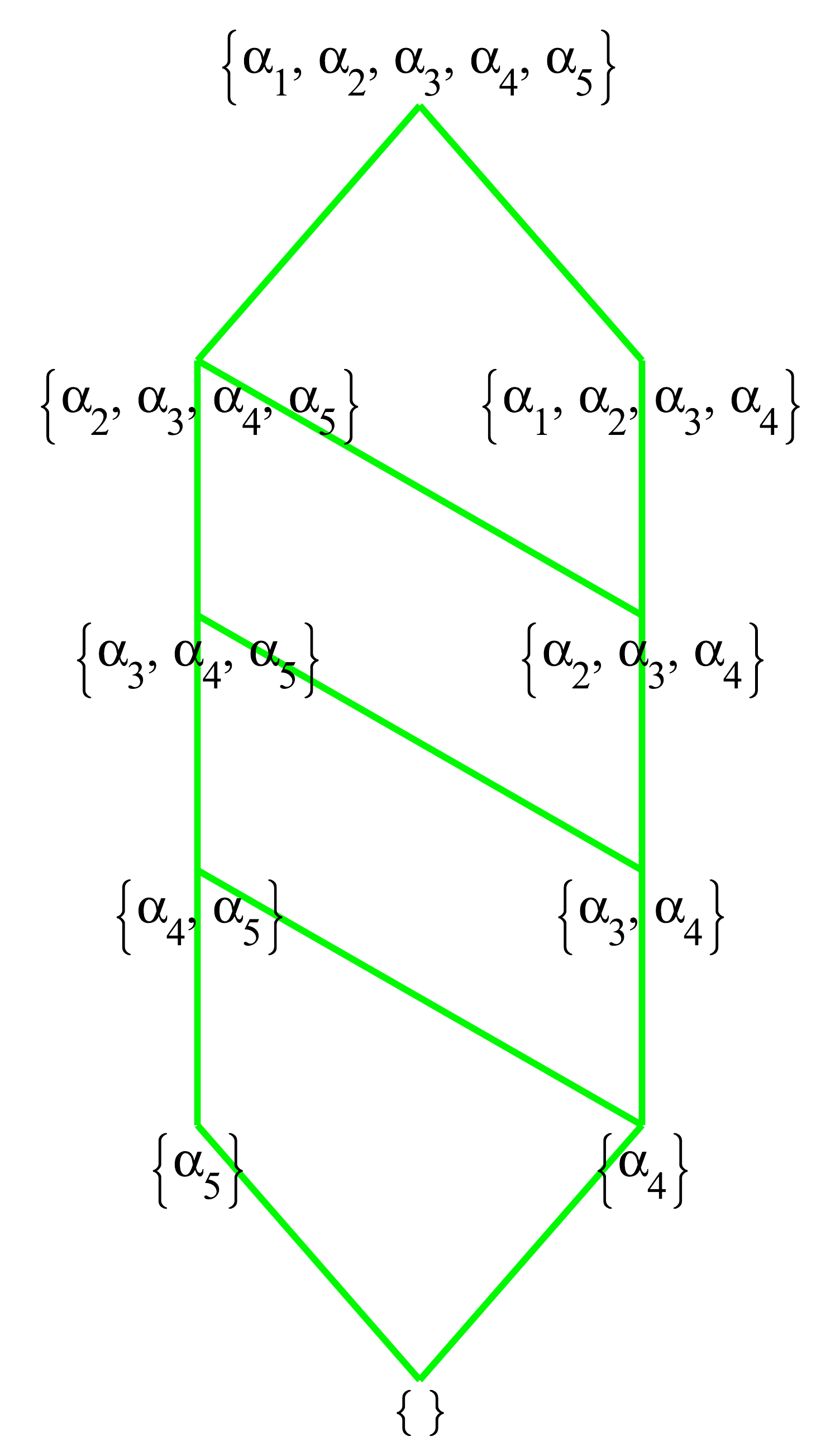}
	      \label{F123}
	  }
	    \subfigure[$J_0=\{\alpha_1,\alpha_5\} $]
	  {
	      \includegraphics[width=2.5in]{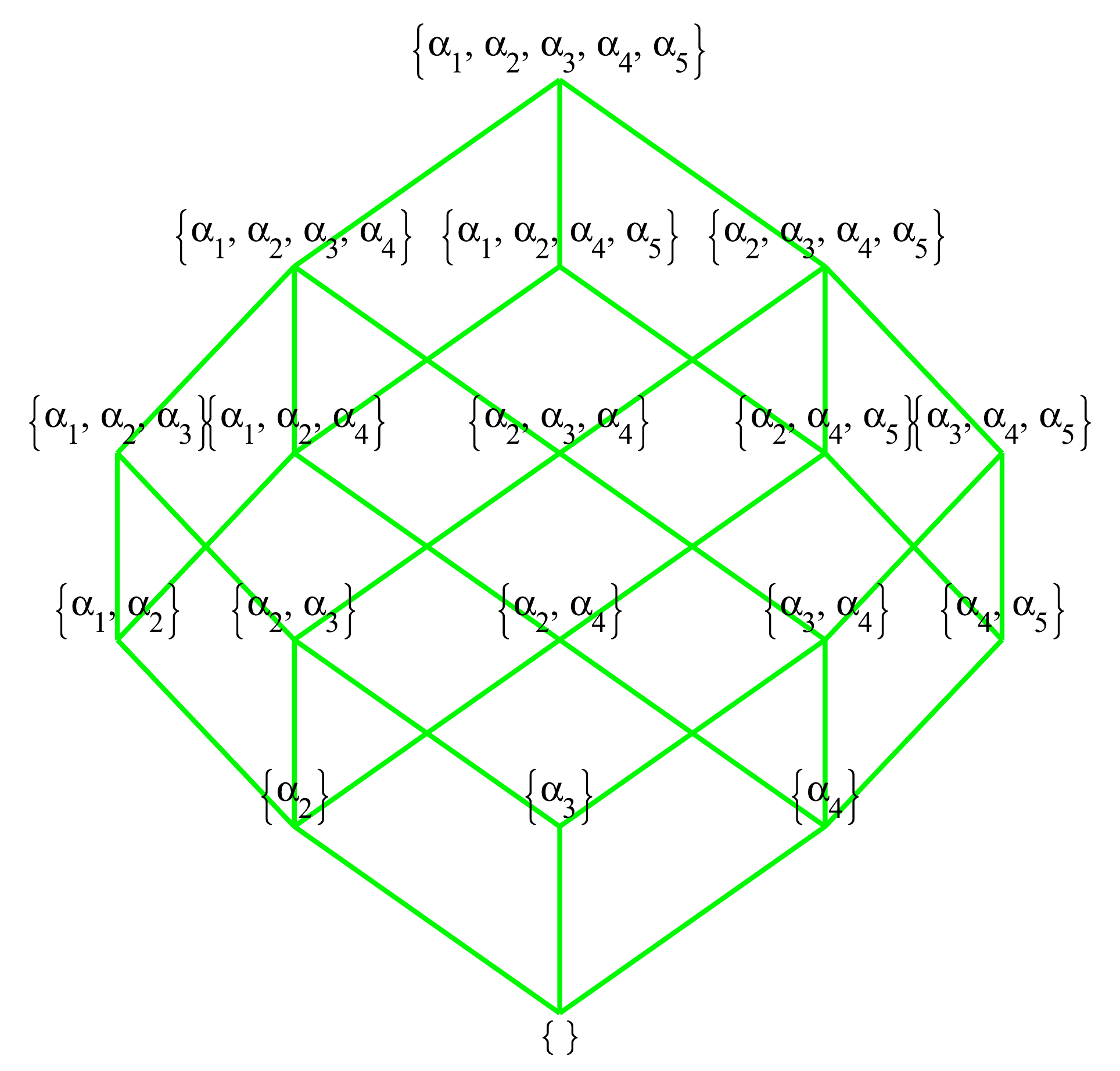}
		      \label{F15}
	  }	
	  \subfigure[$J_0=\{\alpha_1,\alpha_2\} $]
	  {
	      \includegraphics[width=2in]{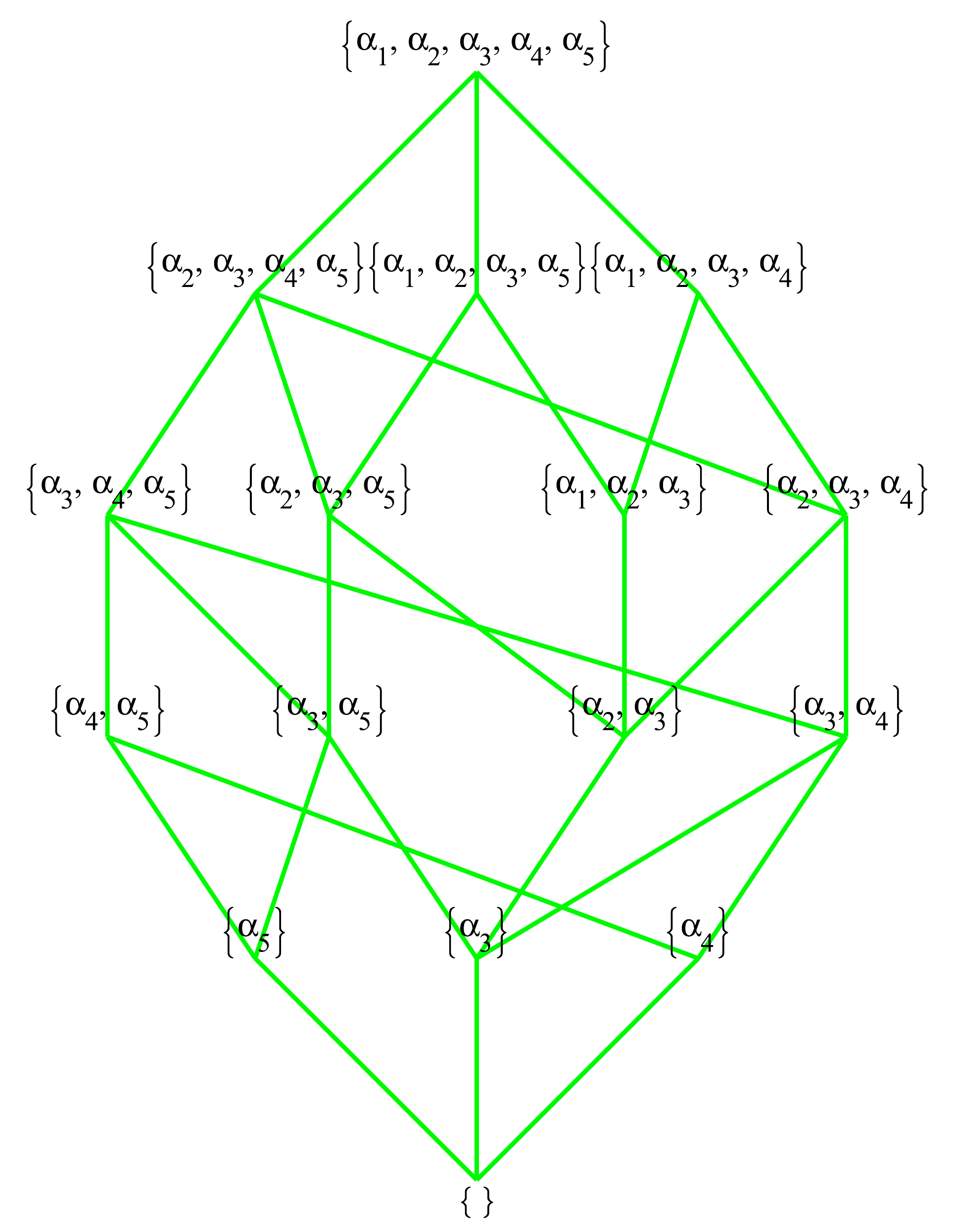}
	      \label{F12}
	  }
	  \caption{Some distributive cross section lattices for $\Delta = \{\alpha_1,...,\alpha_5 \}$.}
	  \label{F:subposets}
	\end{figure}

	\begin{Remark}\label{R:PRD}
	\noindent Let $e_0$ be the minimal nonzero idempotent of $\Lambda$, 
	let $\Delta$ be a set of simple roots and let 
	$J_0= \{ \alpha \in \Delta| \ \sigma_{\alpha} e_0 = e_0 \sigma_{\alpha}\}$ be as before. 
	In \cite{PR88} Putcha and Renner observed that $\LA$ is a distributive lattice if and only if 
	$\Delta - J_0$ is a connected subset of the Dynkin diagram of $G$ 
	(in the sense of the Definition \ref{D:connected}). 
	In Figure \ref{F:subposets} we have depicted some distributive cross section lattices. 
	The example in Figure \ref{F:adjoint23}, however, is not distributive. For more examples see 
	Renner's book \cite{Renner04}.
	\end{Remark}

	We are going to show that a distributive cross section lattice has to 
	be a product of chains. In order to use Stanley's theorem, however, we need to further 
	analyze the structure of $\LA$. We start by strengthening the 
	observation of Putcha and Renner about the distributivity of the cross sections.

	 \begin{Proposition}\label{P:equal}
	 Let $J_0$ and $\Lambda$ be as in Theorem \ref{T:PR88}.
	 The followings are equivalent
	 \begin{enumerate}
	 \item $\LA$ is isomorphic to a sublattice of the Boolean lattice of all subsets of $\Delta$. 
	 \item $\LA$ is distributive.
	 \item $\Delta - J_0$ is connected.
	 \end{enumerate}
	 \end{Proposition}
	 
	 \begin{proof}

	By Remark \ref{R:PRD} the equality of the second and third items is already known. 
	Recall also that $\LA$ is closed under union operation. 
	Therefore, it is enough to prove that $\LA $ is closed under intersections.

	The $J_0$ consists of one or two connected subsets of $\Delta$. 
	We first assume that it has a single connected piece. 
	Without loss of generality we may assume that $|\Delta|=n$ and that 
	$J_0 = \{\alpha_1,...,\alpha_k\}$ for some $k<n$ such that 
	$\sigma_{\alpha_i} \sigma_{ \alpha_{i+1}} \neq \sigma_{\alpha_{i+1}} \sigma_{\alpha_i}$ 
	for $i=1,...,k-1$. Let $I_1 \subseteq \Delta$ be such that no connected component of 
	$I_1$ lie entirely in $J_0$.  
	Therefore, either $I_1 \cap J_0 = \emptyset$, or
	$I_1 \cap J_0 = \{\alpha_i,\alpha_{i+1},...,\alpha_k\}$ and $\alpha_{k+1} \in I_1$.  
	Let $I_2 \subseteq \Delta$ be another such subset.  
	Then, either $I_1 \cap I_2 \cap J_0$ is empty (hence there is nothing to prove), 
	or $I_1 \cap I_2 \cap J_0 = \{ \alpha_j,\alpha_{j+1},...,\alpha_k\}$, for some $j\leq k$, 
	and $\alpha_{k+1} \in I_1 \cap I_2$.
	In other words, no connected component of $I_1 \cap I_2$ lie entirely in $J_0$. 
	Hence, in this case, $\LA$ is closed under intersections.

	Next, we assume that $\Lambda - J_0 = \{\alpha_i:\ k\leq i \leq l\}$ for some 
	$1\leq k < l \leq n$.  Let $I_1$ and $I_2$ be two subsets of $\Delta$ such that 
	no connected components lie entirely in $J_0$.  
	Without loss of generality we may assume that $I_1 \cap I_2 \cap J_0 \neq \emptyset$. 
	Then, either 
	$I_1 \cap I_2 \cap \{\alpha_l,\alpha_{l+1},...,\alpha_n\} = \emptyset$, or 
	$I_1 \cap I_2 \cap \{\alpha_l,\alpha_{l+1},...,\alpha_n\} \neq \emptyset$ and 
	$I_1 \cap I_2 \cap \{\alpha_1,\alpha_{2},...,\alpha_k\} \neq \emptyset$. 
	In both of these cases we may proceed as in the previous paragraph. 
	Therefore, $I_1 \cap I_2$ has no connected component lying entirely in $J_0$. 
	Hence, $\LA$ is closed under intersections.

	\end{proof}

	 \begin{Remark} 
	 A distributive lattice is always modular, but the converse does not need to be true. 
	 A modular lattice is distributive if and only if it does not contain any interval of rank 
	 three which is isomorphic to a ``diamond,'' $M_3 = \{a,b,c,d,e,f\}$ with 
	 $a<b,c,d$ and $b,c,d < e$.  (See Gr\"atzer, \cite{Grat} ).
	 \end{Remark}

	\begin{Theorem}
	Let $\Lambda$ be a modular cross section lattice, and let 
	$J_0 \subseteq \Delta$ be as before. If $|\Delta - J_0| > 1$, 
	then every interval of rank three in $\LA$ is locally rank  symmetric.   
	If $|\Delta - J_0|= 1$, then every interval of rank three in 
	$\LA - \{e_0\}$ is locally rank  symmetric. 
	\end{Theorem}

	\begin{proof}
	The second assertion can be checked from the Fig. 2 of \cite{PR88}. 
	We are going to prove the first assertion.  If $|\Delta - J_0| =2$, 
	then the first assertion can be checked from Fig. 7.1. of \cite{Renner04}. 
	Therefore, we assume that $|\Delta - J_0| > 2$.

	It is enough to prove that this for $J_0=\{\alpha_1,...,\alpha_k\}$, where 
	$\alpha_i \alpha_{i+1} \neq \alpha_{i+1} \alpha_i$ for $i=1,...,k-1$. 
	(The case $J_0=\{\alpha_l,...,\alpha_n\}$ is identical, and the case 
	$J_0=\{\alpha_1,...,\alpha_k,\alpha_l,...,\alpha_n\}$ is similar.)

	It is convenient to identify $\LA$ by its image in $2^{\Delta}$.  
	Let $[U,V] \subseteq 2^{\Delta}$ be an interval of rank 3. 
	Let $V= U \cup \{\alpha_x,\alpha_y,\alpha_z\}$ for some $x<y<z$ from $\{1,...,n\}$.

	If $\{\alpha_x,\alpha_y,\alpha_z\} \cap J_0 \subseteq \{\alpha_x\}$, then 
	$[U,V]$ is isomorphic to the Boolean lattice of subsets of $\{x,y,z\}$, which 
	is locally rank symmetric.  Thus, we may assume that 
	$\{\alpha_x,\alpha_y\}\subseteq \{\alpha_x,\alpha_y,\alpha_z\} \cap J_0$.  
	Then $U \cup \{\alpha_x\}$ can not be of the form $\phi(f)$ for some $f$ covering 
	$\phi(e)=U$. In other words, $U \cup \{\alpha_x\}$ is not contained in the interval $[U,V]$. 
	By the same token, $U\cup \{\alpha_x,\alpha_z\}$ cannot be in $[U,V]$.

	If, in addition, $\alpha_z\in J_0$, then neither $U\cup \{\alpha_y\}$ nor 
	$U\cup \{\alpha_x,\alpha_y\}$ can be in $[U,V]$. In this case, 
	$[U,V] = \{ U, U\cup \{\alpha_z\},U\cup \{\alpha_y,\alpha_z\}, U\cup \{\alpha_x,\alpha_y,\alpha_z\}=V\}$ 
	is a chain, and hence locally rank symmetric.  
	If $\alpha_z \notin J_0$, then 
	$[U,V]= \{U,U\cup \{\alpha_y\}, U\cup \{\alpha_z\},U\cup \{\alpha_x,\alpha_y\}, U\cup \{\alpha_y,\alpha_z\}, 
	U\cup \{\alpha_x,\alpha_y,\alpha_z\}=V\}$ which is also locally rank symmetric. 

	This finishes the proof.
	\end{proof}

	The proof of the theorem gives the proof of the following corollary, also. 
	
	\begin{Corollary}
	Let $\Lambda$ be a distributive (hence modular) cross section lattice of a 
	$J-$irreducible monoid of type $A_{n}$. Let $J_0$ be as before. 
	If $|\Delta 	- J_0| > 1$, then every interval of rank three in 
	$\LA$ is isomorphic to one of the followings.  
	\begin{enumerate}
	\item $2^{\{a,b,c\}}$, the Boolean lattice on three letters.
	\item  $\{a,b,c,d,e,f\}$ where $a<b,\ a<c,\ b<d,\ c<d,\ c<e,\ e<f,\ d<f$.
	\item  $\{a,b,c,d\}$ where $a<b<c<d$.
	\end{enumerate}
	\end{Corollary}

	\begin{Theorem} (F. Regonati, \cite{Regonati}) 
	Let $L$ be a finite modular lattice. The following three conditions are equivalent. 
	\begin{enumerate}
	\item $L$ is locally rank symmetric.
	\item Every interval of $L$ of rank three is rank-symmetric. 
	\item $L$ is a product $P_1 \times P_2 \times \cdots \times P_m$ of $q_i$--primary lattices 
	$P_i$ (including the possibility $q_i=0$, in which case $P_i$ is a chain.
	\end{enumerate}
	\end{Theorem}

	\begin{Corollary}
	Let $\Lambda$ be a distributive (hence modular) cross section lattice of a 
	$J-$irreducible monoid of type $A_{n}$. Let $e_0\in \Lambda$ and 
	$J_0 \subseteq \Delta$ be as before. If $|\Delta - J_0| > 1$, then $\LA$ 
	is locally rank symmetric and is isomorphic to a product of chains. 
	If $|\Delta - J_0| =1$, then $\LA - \{e_0\}$ is locally rank symmetric 
	and is isomorphic to a product of chains. 
	\end{Corollary}

	\begin{Definition}
	Let $\gamma = \gamma_1+\cdots +\gamma_l$ be a partition (of a composition), 
	and let $L_{\gamma}$ be the product 
	$C_{\gamma_1}\times \cdots \times C_{\gamma_l}$ of chains $C_{\gamma_i}$ 
	of length $\gamma_i-1$. We call $\gamma$, the {\em partition type} of $L_{\gamma}$. 
	\end{Definition}

	 It is easy to see that the flag quasi-symmetric function 
	 $F_{L_{\gamma}}$ of $L_{\gamma}=C_{\gamma_1}\times \cdots \times C_{\gamma_l}$ 
	 is equal to the complete homogenous symmetric function 
	 $h_{\gamma} = h_{\gamma_1}\cdots h_{\gamma_l}$ (see Proposition 3.3, \cite{Stan96}).

	It is easy to see that the partition types of the distributive cross section lattices in 
	Figure \ref{F:dortresim} are $(3,2), (4,1), (2,2,1)$ and $(3,1,1)$ in the given order 
	of the figures (a),(b),(c) and (d). Notice that these are partitions of 5. 
	In fact, it is easy to see that, for $n\leq 6$ the number of distinct distributive 
	cross section lattices associated with $GL_n$ is equal to number of partitions of $n-1$. 
	Unfortunately this is not a general fact. The number of distinct distributive cross section 
	lattices associated with $GL_7$ is 10, but the number of partitions of 6 is 11. 
	This reminds us the curious phenomena of ungradedness of the dominance 
	partial order on partitions for $n \geq 6$.

	\begin{Remark}
	Note that the Conjecture \ref{C:chains} gives the possible ``partition types" of the distributive cross section 
	lattices; possible they are of the form $(k+1,n-l+2,1,1...,1)$, $k\geq 0,\ l\geq 0$.  
	\end{Remark}

	\subsubsection{Combinatorially smooth subsets of $W$}

	  An important result of Renner says that, when $M$ is $J$-irreducible monoid, 
	  the poset of idempotents of the (maximal) torus embedding $E(\overline{T})$ is 
	  isomorphic to the face lattice of the polytope $Conv(W\cdot \mu)$, where $\mu$ is 
	  the highest weight associated with the defining irreducible representation of  $M$.

	The set of {\em nonzero minimal elements} in $E(\overline{T})$ 
	$$E_1(\overline{T} ) = \{ e\in E(\overline{T}):\ \dim (Te)=1\}$$ 
	plays a special role in Renner's descent sysyems (see \cite{Renner09}).

	\begin{Definition}
	Let $e, e'\in E_1(\overline{T})$ be two different minimal, nonzero idempotents.  We write 
	$ e < e'$ whenever $eBe' \neq 0$.
	\end{Definition}

	To be consistent with Renner's notation from \cite{Renner09}, in this section, we let $e_1$ denote the  
	minimal nonzero idempotent of $\Lambda$.

	\begin{Theorem}(Renner, \cite{Renner09}) \label{T:minidem}
	The followings are equivalent for $v,w \in W^J$.
	\begin{enumerate}
	\item $e= v e_1 v^{-1} < e'= w e_1 w^{-1}$ in $(E_1,<)$
	\item $w < v $ in $(W^J,<)$ in the Bruhat ordering on $W^J$. 
	\end{enumerate}
	\end{Theorem}

	It is implicit in this theorem that $E_1(\bar{T})$ can be identified with $W^J$. 
	This follows from (\ref{E:ELAMB}) and the fact that $M$ is $J-$irreducible. 
	Recall that $J$ is called combinatorially smooth if the toric variety associated 
	with the polytope $Conv(W\cdot \mu)$ is rationally smooth.

	The reader can extend the following observation to other types: 
	
	\begin{Proposition}	
	Let $\Delta$ be a set of simple roots of type $A_n$, $n\geq 2$.
	Then, a subset $J \subseteq S=\{ \sigma_{\alpha}:\ \alpha \in \Delta\}$ 
	of the simple reflections is combinatorially smooth 
	if and only if $\LA$ (defined by $J_0=J$) is distributive and one of the followings hold: 
	\begin{enumerate}
	\item[a)] $| \Delta - J | \geq 1$, if $J$ contains only one of the end nodes
	of the Dynkin diagram,
	\item[b)] $|\Delta - J| \geq 2$, if $J$ contains both ends of the Dynkin diagram.
	\end{enumerate}
	\end{Proposition}

	\begin{proof}
	See Corollary 3.5 in \cite{Renner09}.
	\end{proof}

	\section{\textbf{Supersolvability}}\label{S:SS}

	In this section we find a necessary and sufficient condition for a cross section lattice of a 
	$J-$irreducible to be supersolvable. 
	We assume for this section that all the cross section lattices belong to one of the following 
	types $A_n,B_n$ or $C_n$.

	\begin{Definition}
	A finite lattice $L$ is called supersolvable if it possesses a maximal chain $\vGam$, 
	called {\em modular chain} or $\mo-chain$, with the property that the sublattice of $L$ 
	generated by $\vGam$ and any other chain of $L$ is distributive. 
	\end{Definition}

	Recall that an element $b\in L$ of a lattice is called right modular if and only if for every $a\in L$ 
	\begin{equation}\label{E:rightmodular}
	c \leq b \Rightarrow c \vee (a \wedge b ) = (c \vee a )\wedge b.
	\end{equation}
	Similarly, $a\in L$ is called left modular if and only if for every 
	$b\in L$, (\ref{E:rightmodular}) holds. 
	An element $a\in L$ is called modular if it is both right and left modular.

	 The following characterization of semismodular supersolvable lattices  is useful.   
	\begin{Lemma} (Stanley, Corollary 2.3, \cite{Stan72})\label{C:USS}
	Let $L$ be a finite upper semimodular lattice, and $\vGam$ be a maximal chain of $L$. 
	Then, $\vGam$ is a modular chain 
	if and only if every element of $\vGam$ is modular. 
	\end{Lemma}

	To prove the Theorem \ref{T:ss} we need the followings, as well:

	\begin{Lemma}\label{L:emptyintersection}
	Let $X$ be the image $\phi(e)\in\phi(\Lambda)$ 
	of an idempotent under the map $\phi: \Lambda \rightarrow 2^{\Delta}$ of the 
	Theorem \ref{T:PR88}.  If $X\cap J_0 = \emptyset$, then $e$ is both right and left modular.
	\end{Lemma}

	\begin{proof}
	
	By abuse of notation, using the map $\phi$ of the Theorem (\ref{T:PR88}), 
	we identify $\Lambda$ by its image.  (Therefore, if $h\leq f$ in $\Lambda$ 
	and $Z=\phi(h),\ Y=\phi(f)$, then we write $Z \leq Y$.) 
	Let $X=\phi(e)$ be as in the hypotheses: $X \cap J_0 = \emptyset$. 
	
	We first show that $X$ is right modular.
	Let $U$ and $V$ be two subsets from $\Lambda$ such that $V \leq X$. 
	Using  Corollary \ref{C:PR88} $U \wedge X = U \cap X$, and there
	exists $H\in \Lambda$ such that $U \cap X = H$. 
	Since $V \vee H = V \cup H$,  $V\vee ( U \wedge X) = V \cup ( U \cap X)$.
	Since $V \subseteq X$ we see that $V \cup ( U \cap X) = (V \cup U) \cap X$.
	On the other hand the right hand side of the last equality is 
	$(V \vee U) \wedge X$.
	Hence for $U,V \in \Lambda$ such that $V \leq X$ 
	the implication in (\ref{E:rightmodular}) holds. 
	In other words, $X$ is a right modular element of $\Lambda$.

	Next, we show that $X$ is left modular. We need to show that for every 
	$V\in \Lambda$ and $U \leq V$, 
	\begin{equation}\label{E:leftwanted}
	U \vee (X \wedge V) = (U \vee X) \wedge V.
	\end{equation}
	 By the Corollary \ref{C:PR88}, the left hand side of the Equation (\ref{E:leftwanted}) 
	 is equal to $U\cup ( X \cap V)$, and $(U\vee X)= (U\cup X)$. 
	 Hence, $(U \vee X) \wedge V \subseteq (U \cup X) \cap V$. 
	 Therefore, it is enough to prove that $ (U \cup X) \cap V \subseteq (U \vee X) \wedge V $.  
	 To this end, let $\alpha \in (U \cup X) \cap V$.  If $\alpha \notin J_0$, then $J_0$ 
	 cannot contain the connected component of $\alpha$, so $\alpha \in (U \vee X) \wedge V$. 
	 Therefore, we may assume that $\alpha \in J_0$.  It follows that $\alpha \in U \subseteq V$.  
	 Therefore, connected component of $\alpha$ in $(U \cup X) \cap V$ is at least as large 
	 as the connected component of $\alpha \in U$. In other words, $J_0$ cannot contain the 
	 connected component of $\alpha$, and hence $\alpha \in (U \vee X) \wedge V$. 
	 This finishes the proof.
	 
	\end{proof}

	\begin{Lemma}\label{L:nonemptyintersection}
	Suppose that a connected component of $J_0$ is either a singleton 
	$\{\alpha_i\}$ , $\alpha_i \in \Delta$, or is of the form 
	$\widetilde{J_0}=\{\alpha_{i_1},\alpha_{i_2},...,\alpha_{i_k}\} \subseteq J_0,\ k>1$ such that
	\begin{enumerate}
	\item if $\alpha_{i_k} \in \widetilde{J_0}$, then there exists $\alpha_{i_l} \in \widetilde{J_0}$ 
	such that $\sigma_{\alpha_{i_k}} \sigma_{\alpha_{i_l}} \neq \sigma_{\alpha_{i_l}} 
	\sigma_{\alpha_{i_k}}$, and
	\item there exists an end-node $\alpha_{i}$ of $\Delta$ contained in $\widetilde{J_0}$. 
	\end{enumerate}
	Let $X=\phi(e) \in\phi(\Lambda)$ be such that $\Delta - J_0 \subseteq X$. 
	Then,  $e$ is both right and left modular.
	\end{Lemma}

	\begin{proof}
	We proceed as in the proof of the Lemma \ref{L:emptyintersection}, 
	using the map $\phi$ of the Theorem (\ref{T:PR88}) to identify 
	$\Lambda$ by its image. 
	Let $X$ be as in the hypotheses, so that $\Delta - J_0 \subseteq X$. 
	Notice if we can show that for every $U\in \Lambda$, $U\wedge X = U \cap X$, 
	then we are done by the proof of the Lemma \ref{L:emptyintersection}.

	To this end, we assume that there exists $U\in \Lambda$ such that 
	$U \wedge X \neq U \cap X$. Thus, by Corollary \ref{C:PR88},  there must exist 
	$\alpha \in J_0 \cap U \cap X$ with a connected neighborhood (in $U \cap X$) 
	which lies completely in $J_0$.  
	Since $\Delta - J_0 \subseteq X$, if $\{\alpha \}\in J_0$ is a singleton subset of 
	$J_0$, any connected component of $U$ which contains $\alpha$ has to intersect 
	$X$ with more than one element. In other words, the connected component in 
	$U\cap X$ of $\alpha$ cannot lie in $J_0$ completely.  
	If $\alpha$ lies in a component $\widetilde{J_0}$ of $J_0$ which contains an 
	end-node, then we proceed as in the proof the Proposition \ref{P:equal}. 
	Suppose that $\widetilde{J_0}= \{\alpha_{j_1},\alpha_{j_2},...,\alpha_{j_k}\}$ is 
	the connected subset $J_0$ containing $\alpha$, and $\alpha_{j_k} \in \widetilde{J_0}$ 
	is an end-node (of $\Delta$).

	Then, $U\cap J_0=\{ \alpha_{j_1},\alpha_{j_2},...,\alpha_{j_l} \}$, for some 
	$l \leq k$ such that $\sigma_{\alpha_{j_m}} \sigma_{\alpha_{j_{m+1}}} \neq 
	\sigma_{\alpha_{j_{m+1}}} \sigma_{\alpha_{j_m}}$ for $m=1,...,l-1$.  
	Furthermore, there exists $\alpha' \in U- J_0$ such that 
	$\sigma_{\alpha'} \sigma_{\alpha_{j_1}} \neq \sigma_{\alpha_{j_1}} \sigma_{\alpha'}$.

	Since $\Delta - J_0 \subseteq X$, $\alpha' \in X$. 
	Therefore, $\{ \alpha' \} \cup U\cap J_0 \subseteq X$, and hence, no connected 
	component of $\alpha$ in $U \cap X$ can lie completely in $J_0$.  
	In other words, $U\cap X = U\wedge X$. 
	The rest of the proof goes as in the previous Lemma \ref{L:emptyintersection}. 
	\end{proof}

	\noindent {\em Proof of Theorem \ref{T:ss}.} ($\Leftarrow$) 
	We explicitly construct an $\mo$-chain.  
	Once again, we identify $\Lambda$ with its image in $2^{\Delta}$.
	Let $J_0=\{\alpha_{i_1},\alpha_{i_2},...,\alpha_{i_k}\} \subseteq \Delta$ be as in the 
	hypotheses of the Theorem \ref{T:ss}, and let $\{\alpha_{j_1},...,\alpha_{j_m}\}$ 
	be the complement of $J_0$ in $\Delta$. 
	We assume that $i_1<i_2 < \cdots < i_k$, and $j_1<j_2<\cdots < j_m$.

	By Corollary \ref{C:PR88}, any entry $U$ of the chain 
	$$ 
	\emptyset \subset \{\alpha_{j_1}\} \subset \{\alpha_{j_1},\alpha_{j_2}\} 
	\subset \cdots \subset \{\alpha_{j_1},...,\alpha_{j_m}\}
	$$
	is an element of $\Lambda$, and furthermore, by Lemma \ref{L:emptyintersection}, 
	$U$ is both left and right modular.

	Let $I=\{i_1,...,i_{k'}\} \subseteq [k]$ be the set of indices of elements of 
	$J_0$ which are less than $j_1$.  
	The set $I$ might be empty.  If not, by the hypotheses of the Theorem, 
	$\{i_1,...,i_{k'}\}=\{1,...,k'\}$.
	Then, it is easy to check that the entries of the chain
	\begin{equation} \label{E:orta}
	\{\alpha_{i_{k'}} \} \cup (\Delta - J_0)  \subset \cdots \subset
	\{\alpha_{i_1},\alpha_{i_2}, \ldots , \alpha_{i_{k'}} \} \cup (\Delta - J_0)
	\end{equation}
	as well as the entries of the chain 
	\begin{equation}\label{E:son}
	U_{i_{k'+1}} \subset U_{i_{k'+2}} \subset \cdots \subset U_{i_k} = \Delta, 
	\end{equation}
	where
	\begin{align*}
	U_{i_{k'+1}} &= \{\alpha_{i_1},...,\alpha_{i_{k'}} \}  \cup (\Delta - J_0) \cup \{ \alpha_{i_{k'+1}}\} \\
	U_{i_{k'+2}} &= \{\alpha_{i_1},...,\alpha_{i_{k'}} \}  \cup (\Delta - J_0) \cup \{ \alpha_{i_{k'+1}}, \alpha_{i_{k'+2}}  \} \\
	& \vdots  \\
	U_{i_{k}} &= \{\alpha_{i_1},...,\alpha_{i_{k'}}   \}  \cup (\Delta - J_0) \cup \{  \alpha_{i_{k'+1}},\alpha_{i_{k'+2}}, \ldots, \alpha_{i_k}\}
	\end{align*}
	are in $\Lambda$.

	By Lemma \ref{L:nonemptyintersection}, any entry of the chain 
	(\ref{E:orta}) and any entry of (\ref{E:son}) is both left and right modular. 
	Therefore, we have found a maximal chain $\vGam$ whose entries are 
	both left and right modular.

	($\Rightarrow$) Assume $\Lambda$ is supersolvable and that  there exists a 
	connected component $\widetilde{J_0} \subseteq J_0$ which does not contain 
	an end-node of $\Delta$, and $|\widetilde{J_0}| > 1$, say 
	$\widetilde{J_0} = \{ \alpha_{i_1},...,\alpha_{i_k}\}$, $k>1$.  
	Then, there exist $\alpha, \alpha' \in \Delta - J_0$ such that 
	$\sigma_{\alpha} \sigma_{\alpha_{i_1}} \neq \sigma_{\alpha_{i_1}} \sigma_{\alpha}$, and  
	$\sigma_{\alpha'} \sigma_{\alpha_{i_k}} \neq \sigma_{\alpha_{i_k}} \sigma_{\alpha'}$. 
	Let $\vGam$ be a modular chain for $\Lambda$. 
	Then,  there exist an entry $\vGam_t$ of $\vGam$ such that there exists 
	$\alpha_{i_{m}} \in \widetilde{J_0} -  \vGam_t$ for some $1 \leq m < k$. 
	Without loss of generality we may assume that $\alpha_{i_n} \in \vGam_t$ 
	whenever $m < n \leq k$.
	Let $C$ be the subset $\{ \alpha_{i_{m+2}},...,\alpha_k,\alpha'\} \subseteq \vGam_t$. 
	Let $A$ be the set $\{ \alpha, \alpha_{i_1},...,\alpha_{i_m},\alpha_{i_{m+1}}\}$.  
	Clearly, $A$ and $C$ are in $\Lambda$. Since $\alpha_{i_{m+1}} \in A \cap \vGam_t$ 
	is isolated in $J_0$, 
	$A \wedge \vGam_t$ cannot contain $\alpha_{i_{m+1}}$. 
	Therefore, $\alpha_{i_{m+1}} \notin C \vee (A \wedge \vGam_t)$. 
	However, it is easy to check that $\alpha_{i_{m+1}} \in (C \vee A)\wedge \vGam_t$. 
	In other words, 
	$$
	C \vee ( A  \wedge \vGam_t) \neq (C \vee A) \wedge  \vGam_t.
	$$
	Therefore, the element $\vGam_t$ of the modular chain $\vGam$ is not (right) modular, 
	which contradicts with the Lemma \ref{C:USS}. 
	Therefore, we must have $|\widetilde{J_0}| \leq 1$. This finishes the proof.

	\section{\textbf{Characteristic polynomials}}\label{S:CP}

	Recall that the M\"obius function of a poset $P$ is the unique function 
	$\mu: P\times P \rightarrow \N$ satisfying 
	\begin{enumerate}
	\item $\mu(x,x) =1$ for every $x\in P$,
	\item $\mu(x,y)=0$ whenever $x \nleq y$,
	\item $\mu(x,y) =- \sum_{x\leq z < y} \mu(x,z)$ for all $x<y$ in $P$.
	\end{enumerate}
	The {\em characteristic polynomial} $p(\alpha,P)$ (also known as 
	{\em Birkhoff polynomial}) of a finite graded poset $P$ of rank $n$ is 
	\begin{equation}
	p(x,P) = \sum_{x\in P} \mu(\hat{0},x) x^{rk(\hat{1})-rk(\hat{x})}.
	\end{equation}
	A particularly nice survey about characteristic polynomials is written by 
	B. Sagan and can be found at \cite{Sag99}.

	\begin{Conjecture}
	Let $\Lambda$ be the cross section lattice of a $J-$irreducible monoid $M$. Then
	\begin{equation}
	p(x,\LA)=  x^{|J_0|} (x-1)^{n-|J_0|}.
	\end{equation}
	\end{Conjecture}

	In this section we prove the following following special case of the above conjecture. 
	\begin{Theorem}\label{charpoly}
	Let $n$ be the rank of the supersolvable cross section lattice $\LA$. 
	Then the characteristic polynomial of $\vGam$ is 	
	\begin{equation}
	p(x,\LA)=  x^{|J_0|} (x-1)^{n-|J_0|}.
	\end{equation}
	\end{Theorem}

	\begin{proof}

	Recall Theorem \ref{T:Stanley72}: 
	If  $L$ be a semimodular supersolvable lattice and  
	$\hat{0} = x_0 < x_1 <\cdots < x_n = \hat{1}$ is an $\mathcal{M}$-chain,
	then
	\begin{equation}
	p(x,L) = (x - a_1)(x-a_2) \cdots (x - a_n),
	\end{equation}
	where $a_i$ is the number of atoms $u\in L$ such that $u \leq x_i$ and $u \nleq  x_{i-1}$.

	In the proof of Theorem  \ref{T:ss} we found 
	an $\mathcal{M}$-chain $D$.  It is easy to check that $\Delta - J_0$ 
	is the set of atoms of $\LA$. Furthermore, initial part 
	$x_0=\hat{0} < x_1 < \cdots < x_k$ of $D$ is given by the subsets 
	$x_s = \{j_1,...,j_s\}_< \subseteq \Delta - J_0$. 
	Therefore, $a_s = 1$ for $s=1,...,|\Delta - J_0|$. 
	Since $x_s \subseteq x_r$ for $s<r$, it follows that $a_r = 0$. 
	Therefore, 
	\begin{equation}
	p(x,\LA) = \prod_{s=1}^{|\Delta - J_0|} (x - 1) 
	\prod_{r=|\Delta - J_0|+1}^{|\Delta|} x =x^{|J_0|} (x-1)^{n-|J_0|}.
	\end{equation}
	This finishes the proof. 

	\end{proof}

	The monoid in the following example is interesting in its own right. 
	
	\begin{Example}
	 Let $\Lambda$ be the cross section lattice of the monoid obtained from 
	 $SL_5$ by using the adjoint representation $\rho(g) = \text{Ad} (g)$. 
	 Then, one can show that 
	 $J_0 = \{\alpha_2,\alpha_3\} \subseteq \Delta = \{\alpha_1,\alpha_2,\alpha_3,\alpha_4\}$. 
	The Hasse diagram of $\LA$ is as in Figure \ref{F:adjoint23} 
	(see also Section 7.4.1. of \cite{Renner04}).
	By Theorem \ref{T:ss}, $\LA$ is not supersolvable, however, 
	a straightforward calculations shows that the characteristic polynomial of 
	$\LA$ is equal to $x^2(x-1)^2$.

	\begin{figure}[ht]
	\begin{center}
	\includegraphics[ scale=0.3]{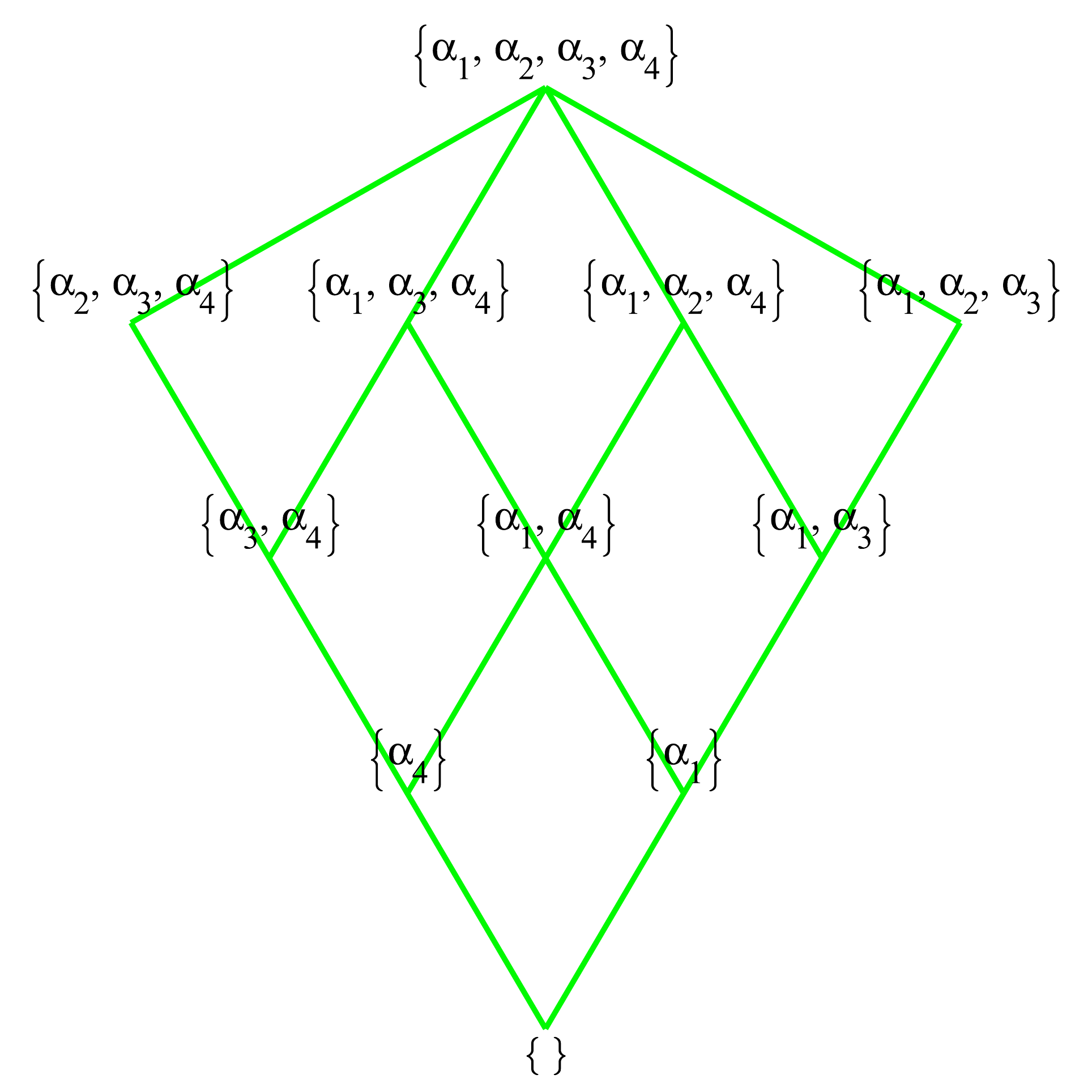}
	\caption{The cross section lattice of an adjoint representation.}
	\label{F:adjoint23}
	\end{center}
	\end{figure}

	\end{Example}

	\section{\textbf{Future directions}}

	In this section, we briefly report some of our progress as a continuation of the 
	results of this article. However, we take the liberty of not introducing 
	the further notation.

	In \cite{McNamara03}, McNamara shows that a lattice $L$ is supersolvable 
	if and only if there exists a 0-Hecke algebra action on the maximal chains of the lattice $L$. 
	He furthermore shows that under a suitable analogue of the Frobenius characteristic map, 
	the representations on the maximal chains might be identified with the flag quasi 
	symmetric functions studied by Stanley (see Section \ref{S:posetterminology}, 
	and \cite{Stan96}). We will write our progress on the representations of the 
	0-Hecke algebra on $\LA$ in a forthcoming article. However, we would like to 
	mention that, in general, flag quasi symmetric functions of Stanley are nice for a cross section
	lattice of a $J-$irreducible monoid. 
	To give the flavor of the results in that direction we state without a proof the following observation:

	 Let $\langle, \rangle$ be the nondegenerate inner product on the space of 
	 quasi symmetric functions such that the set $\{F_{I,n}:\ I\subseteq \{1,...,n-1\}\}$ 
	 of fundamental quasi symmetric functions forms an orthonormal basis.

	\begin{Theorem}
	Let $\LA$ be a rank $n$ cross section a $J-$irreducible monoid. 
	Let $J_0\subseteq \Delta$ be as in Theorem \ref{T:PR88}, and let $n=|\Delta|$.Then, 
	$$
	\langle F_{\LA}, F_{\{1\},n} \rangle = |n-J_0|.
	$$ 
	\end{Theorem}

	The NilCoxeter algebra of a Weyl group plays an important role in the theory of 
	symmetric functions, especially in type $A_n$. 
	We can show that, like a 0-Hecke algebra, the nilCoxeter algeba acts on the 
	maximal chains of $\LA$, also. Thereby, we apply the work of 
	S. Fomin and C. Greene, \cite{FominGreene}. 
	This enables us to generalize Stanley symmetric functions, as well as stable 
	Grothendieck polynomials (via 0-Hecke algebra action) to the setting of the cross section lattices. 
	We will report on these considerations in future papers.

	\subsection{\textbf{Variation of the theme}}

	As a result of Theorem \ref{T:PR88}, the following definition of 
	a ``combinatorial cross section lattice'' is appropriate:

	\begin{Definition}\label{D:combcross}
	Let $\G$ be a graph. By abuse of notation we use $\G$ 
	to denote the set of vertices as well (thus $2^G$ is the 
	set of all subsets of the vertex set $\G$).
	Let $J_0 \subseteq \G$ be a subset of the vertex set.  
	The combinatorial cross section lattice 
	$\Lambda= \Lambda(\G,J_0) \subseteq 2^\G$ associated with the pair $(\G,J_0)$ 
	consists of those subsets $U\in 2^\G$
	having no connected component contained entirely in $J_0$ as an induced 
	subgraph.  The partial ordering on $\Lambda$ is the set inclusion. 
	We consider the empty set as an element of $\Lambda$. 
	\end{Definition}

	\begin{Remark}
	It is clear that if $U$ and $V$ are from 
	a combinatorial cross section poset $\Lambda$, then 
	the join $U \vee V \in \Lambda$ exists and equal to the union $U \cup V$.
	If $U\cap V$ is an element of $\Lambda$, then it is equal to the meet $U \wedge V$,
	otherwise the meet is equal to $\emptyset \in \Lambda$.
	\end{Remark}

	Obviously, when $\G$ is the graph $\G = \{\alpha_1,\ldots, \alpha_n\}$ 
	with the edge set $E=\{ e_1,\ldots, e_{n-1} \}$, where $e_i$ connects 
	$\alpha_i$ and $\alpha_{i+1}$ an associated combinatorial cross section lattice 
	is equal to a cross section lattice of a $J$-irreducible monoid of type 
	$A_n$, considered in the manuscript.

	We are planning to investigate the case of an arbitrary graph 
	in a future paper. However, let us mention here briefly the case when $\G$ 
	is a circuit. In other words, the (Coxeter) graph of an affine root system 
	of type $\tilde{A}_{n}$.
	Its vertex set is $\G= \Delta_{af} = \Delta \cup \{ \alpha_0 \}$,
	where $\Delta$ is a set of simple roots of type $A_{n}$, and $\alpha_0$ is connected
	to both of the end nodes of $\Delta$.

	\begin{Theorem}\label{T:circulartoline}
	Let $\G$ be a circuit, and let $J_0 \subseteq \G$. 
	\begin{enumerate}
	\item If $\G - J_0$
	contains (at least) two adjacent vertices, then $\Lambda(\G,J_0)$ is 
	isomorphic to a cross section lattice of a $J$-irreducible monoid of type $A$.  
	\item $\G$ is supersolvable if and only if each connected component of 
	$J_0$ is a singleton.
	\end{enumerate} 
	\end{Theorem}

	\begin{proof}
	\begin{enumerate}
	\item Suppose that the vertex set of $\G$ is $\{\beta_0,\beta_1,\ldots, \beta_n \}$
	and its edge set is $\{ e_0,\ldots, e_n\}$ where $e_i$ connects $\beta_i$
	and $\beta_{i+1}$, $i=0,\ldots, n-1$, and $e_n$ connects $\beta_n$ and $\beta_0$.   
	Without loss of generality we may assume that $\beta_0$ and $\beta_1$ 
	are the two vertices which are not contained in $\G - J_0$. 
	
	Let $\Delta $ be a root system (of type $A_{n+1}$) 
	with the set of simple roots $\{ \alpha_1,\ldots, \alpha_{n+1}\}$ such that 
	there exists an edge between $\alpha_i $ and $\alpha_{i+1}$ for $i=1,\ldots, n$. 
	Define $$J_0'  =\{  \alpha_i \in \Delta:\  \beta_i \in J_0  \}$$

	For $U \in \Lambda( \G, J_0)$ define $\Phi (U)$ by 
	\begin{equation*}
	\alpha_i \in \Phi(U)\ \text{if}\ 
	\begin{cases}
	i \neq n+1\ \text{and}\ \beta_i \in U,\\
	i= n+1\ \text{and}\ \beta_0 \in U. 
	\end{cases}
	\end{equation*}
	
	Let $\Lambda$ be the cross section lattice of a $J$-irreducible 
	monoid on $\Delta$ determined by $J_0'$ as in the Theorem \ref{T:PR88}. 
	Then, it is easy to check that $U \mapsto \Phi(U)$ is an isomorphism between
	$\Lambda( \G, J_0)$ and $\Lambda$.

	\item Similar to the proof of the Theorem \ref{T:ss}.
	
	\end{enumerate}
		
	\end{proof}
	
	\begin{Remark}
	It is easy to see that the proof of the first part of the Theorem \ref{T:circulartoline} can be reversed 
	to show that if $J_0'\subseteq \Delta$ does not contain any of the end nodes of the 
	diagram of $\Delta$, then $\Lambda$ is isomorphic to a combinatorial cross section lattice 
	$\Lambda(\G,J_0)$ for some $J_0 \subset \G$, and $\G$ is a circuit. 
	\end{Remark}

\bibliographystyle{amsalpha}
 
\end{document}